\documentclass{amsart}
\usepackage{amsmath,a4wide,amssymb,xcolor,pict2e}
\usepackage{graphics,bm}
\usepackage{hyperref}
\usepackage[all]{xy}
\newcommand{\subparallel}{\,|\mbox{\hskip-1pt}\big|\,}
\def\defterm#1{{\color{purple}\em #1}}
\newcommand{\defeq}{:=} 
\newcommand{\Ra}{\Rightarrow}
\newcommand{\IN}{\mathbb N}
\newcommand{\IR}{\mathbb R}
\newcommand{\IZ}{\mathbb Z}

\newcommand{\IF}{\mathbb F}
\newcommand{\F}{\mathcal F}

\newcommand{\w}{\omega}

\newcommand\Af{\mathsf L}

\def\vecbold#1{\vec{\boldsymbol {#1}}}

\def\Aline#1#2{\overline{#1\phantom{\mbox{\footnotesize$|$}}#2}}

\title{Linear Geometry: flats, ranks, regularity, parallelity}
\author{Taras Banakh, Ivan Hetman, Alex Ravsky, Vlad Pshyk}
\address{T.~Banakh: Ivan Franko National University of Lviv (Ukraine) \newline\indent and Jan Kochanowski University in Kielce (Poland)}
\email{t.o.banakh@gmail.com}
\address{I.~Hetman: Lviv (Ukraine)}
\email{vespertilion@gmail.com}
\address{A.~Ravsky: Ya. Pidstryhach Institute for Applied Problems of Mechanics and Matematics, Lviv}
\email{vitaravski@gmail.com}
\address{V.~Pshyk:  Ivan Franko National University of Lviv (Ukraine)}
\email{pshik.vlad8@gmail.com}

\newtheorem{theorem}{Theorem}[section]
\numberwithin{theorem}{section}
\newtheorem{corollary}[theorem]{Corollary}
\newtheorem{proposition}[theorem]{Proposition}
\newtheorem{lemma}[theorem]{Lemma}

\newtheorem{problem}[theorem]{Problem}

\newtheorem{question}[theorem]{Question}

\theoremstyle{definition}

\newtheorem{remark}[theorem]{Remark}
\newtheorem{definition}[theorem]{Definition}
\newtheorem{example}[theorem]{Example}

\begin{document}

\begin{abstract} Linear Geometry describes geometric properties that depend on the fundamental notion of a line. In this paper we survey basic notions and results of Linear Geomery that depend on the flat hulls: flats, exchange, rank, regularity, modularity, and parallelity.
\end{abstract} 
\maketitle

\section{Introduction}
This is a first part of a series of surveys in Linear Geometry, an area of Geometry studying geometric properties that depend on lines only. Linear Geometry has been developed by many mathematicians since times of antiquity. It was formed as a separate area of Geometry in XX-th century.
We present the main ideas and results of Linear Geometry using Bourbaki approach of mathematics which studies mathematical structures. The fundamental mathematical structure of Linear Geometry is that of a liner.

\begin{definition} A \defterm{liner} is a set $X$ whose elements are called \defterm{points}, endowed with a family $\mathcal L$ of subsets of $X$, called \defterm{lines}, that satisfy two axioms:
\begin{itemize}
\item[(L1)] any distinct points belong to a unique line;
\item[(L2)] every line contains at least two distinct points.
\end{itemize}
The family of lines $\mathcal L$ is called the \defterm{line structure} of the liner $X$. 
\end{definition}

This seemingly poor geometric structure allows to develop a rich geometric theory, full of fascinating theorems and intriguing open problems. A self-contained presentation of this theory can be found in the book \cite{Banakh} of the first author. In this paper we survey main notions and results of Chapters 1--6 and 25 of the book \cite{Banakh}. In those chapters we concider concepts that are expressible in via flat hulls: ranks, modularity and parallelity.




\section{Flat sets in liners}

\begin{definition} A subset $A$ of a liner $X$ is called \defterm{flat} if for any distinct points $x,y\in A$, the unique line $\Aline xy$ containing the points $x,y\in A$ is contained in the set $A$. 
\end{definition}

This definition implies that every subset $A\subseteq X$ of cardinality $|A|\le 1$ is flat. 

The uniqueness of a line passing through two distinct points of a liner implies that every line in a liner is a flat set. 

\begin{definition} For a subset $A$ of a liner $X$, the smallest flat set containing $A$ is called the \defterm{flat hull} of $A$ and is denoted by $\overline A$. 
\end{definition}

The flat hull of the empty set is empty; the flat hull of a singleton is this singleton, and the flat hull of any doubleton $\{x,y\}$ is the line $\Aline xy$ passing through the points $x,y$.

\begin{proposition}\label{p:flat-hull} For any subset $A$ of a liner $X$, the flat hull $\overline A$ is equal to the union $\bigcup_{n\in\w}A_n$ of the sets defined by the resursive formula: $A_0\defeq A$ and $A_{n+1}=\bigcup_{x,y\in A_n}\Aline xy$ for $n\in\w$.
\end{proposition}

Here we assume that $\Aline xy=\{x\}=\{y\}$ if $x=y$.

Proposition~\ref{p:flat-hull} implies the following description of flat hulls in liners.

\begin{corollary}[{\cite[1.5.2]{Banakh}}] For any set $A$ in a liner $X$, its flat hull $\overline A$ equals $\bigcup_{F\in[A]^{<\w}}\overline F$.
\end{corollary}

Here $[A]^{<\omega}$ stands for the family of all finite subsets of a set $A$.
For a set $X$ we denote by $|X|$ the cardinality of $X$.

\begin{definition} The \defterm{rank} $\|A\|$ of a subset $A$ of a liner $X$ is the smallest cardinality $|B|$ of a subset $B\subseteq X$ whose flat hull $\overline B$ contains the set $A$. It is clear that the rank is monotone in the sense that $\|A\|\le\|B\|$ for any sets $A\subseteq B$ in a liner.
\end{definition}

\begin{definition} A flat subset of rank 3 in a liner $X$ is called a \defterm{plane} in $X$. 
\end{definition}

\begin{definition}A subset $A$ of a liner $X$ is called a \defterm{hyperplane} in $X$ if $A\ne X$ is a flat set such that $\overline{A\cup\{x\}}=X$ for all $x\in X\setminus A$.
\end{definition}

\begin{definition}  Every subset $A$ of a liner $(X,\mathcal L)$ carries the induced line structure $$\mathcal L{\restriction}_A\defeq\{L\cap A:L\in\mathcal L\;\wedge\;|L\cap A|\ge 2\}.$$  The liner $(A,\mathcal L{\restriction}_A)$ is called a \defterm{subliner} of the liner $(X,\mathcal L)$.
\end{definition}

\section{Exchange Properties in liners}\label{s:exchange}

\begin{definition}  A liner $X$ is defined to have the \defterm{Exchange Property} if for every flat $A\subseteq X$ and points $x\in X\setminus A$, $y\in \overline{A\cup\{x\}}\setminus  A$ we have $x\in\overline{A\cup\{y\}}$. \end{definition}

The Exchange Property first appeared in the Steinitz Exchange Lemma, used by Steinitz for introducing the dimension of vector spaces.
 The Exchange Property can be quantified as follows.

\begin{definition}\label{d:k-EP} A liner $X$ has a \defterm{$\kappa$-Exchange Property} for a cardinal $\kappa$ if for any set $A\subseteq X$ of cardinality $|A|<\kappa$ and any points $x\in X$ and $y\in\overline{A\cup\{x\}}\setminus\overline A$, we have $x\in\overline{A\cup \{y\}}$.
\end{definition}

It is clear that a liner has the Exchange Property if and only if it has the $\kappa$-Exchange Property for every cardinal $\kappa$.  Every liner has the $2$-Exchange Property. 

\begin{proposition}[{\cite[2.1.6]{Banakh}}]\label{p:fh-EP<=>wEP} A liner has the Exchange Property if and only if it has the $n$-Exchange Property for every $n\in\w$.
\end{proposition}

\begin{proposition}[{\cite[2.2.3]{Banakh}}]\label{c:rank-EP} If a liner $X$ has the $\kappa$-Exchange Property for some cardinal $\kappa$, then for every subset $A\subseteq X$ of rank $\|A\|\le \kappa$, there exists a set $A'\subseteq A$ of cardinality $|A'|=\|A\|$ such that $A\subseteq\overline {A'}$.
\end{proposition}

\begin{example}[Terence Tao\footnote{See ({\tt mathoverflow.net/a/450472/61536}).}]\label{ex:Tao} There exists a liner  that does not have the $3$-Exchange Property.
\end{example}

\begin{proof}[Proof] Let $\IF_3=\{-1,0,1\}$ be a 3-element field, $n,m$ be positive integers and $f:\IF_3^n\to \IF_3^m$ be an even function with $f(\{0\}^n)=\{0\}^m$.  Endow the set $X=\IF_3^n\times \IF_3^m$ with the family of lines 
$$
\begin{aligned}
\mathcal L&=\big\{\{(x,x'),(y,y'),(z,z')\}\in [X]^3:x+y+z=0\;\wedge\;x'+y'+z'=f(x-z)\big\}.
\end{aligned}
$$
Here $[X]^3=\{(x,y,z)\in X^3:|\{x,y,z\}|=3\}$.
It can be shown that  $X$ is a liner. Now we select a special even function $f$ for which this liner does not have $3$-Exchange Property. For $n=2$ and $m=1$, consider the even function $f:F_3^2\to \{0,1\}\subseteq F_3$ such that $f^{-1}(1)=\{(0,1),(0,-1)\}$. For such a function $f$, the set $L=\{(x,0,0):x\in F_3\}$ is a line in $(X,\Af)$ such that for the point $\vecbold{u}=(0,1,0)$ the flat hull $\overline{L\cup\{\vecbold{u}\}}$ coincides with $X$. On the other hand, for the point $\vecbold{v}=(0,0,1)\in\overline{L\cup\{\vecbold{u}\}}\setminus L$ the flat hull $\overline{L\cup\{\vecbold{v}\}}$ coincides with the set $F_3\times\{0\}\times F_3$, which is strictly smaller than $\overline{L\cup\{\vecbold{u}\}}=X$. This means that $X$ does not have the $3$-Exchange Property.

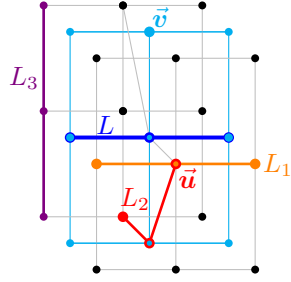
\begin{picture}(100,120)(-180,-60)

\put(-30,40){\color{cyan}\line(1,0){60}}
\put(-30,-40){\color{cyan}\line(1,0){60}}

\put(0,-40){\color{cyan}\line(0,1){80}}
\put(20,-30){\color{lightgray}\line(0,1){80}}
\put(30,-40){\color{cyan}\line(0,1){80}}
\put(40,-50){\color{lightgray}\line(0,1){80}}
\put(-20,-50){\color{lightgray}\line(0,1){80}}
\put(-30,-40){\color{cyan}\line(0,1){80}}
{\linethickness{1pt}
\put(-40,-30){\color{violet}\line(0,1){80}}
}
\put(-53,20){\color{violet}$L_3$}
\put(-10,-30){\color{lightgray}\line(0,1){80}}
\put(10,-50){\color{lightgray}\line(0,1){80}}

\put(-40,10){\color{lightgray}\line(1,0){60}}
{\linethickness{1.5pt}
\put(-30,0){\color{blue}\line(1,0){60}}
}
\put(-20,-10){\color{lightgray}\line(1,0){60}}

\put(-40,50){\color{lightgray}\line(1,0){60}}
\put(-30,40){\color{cyan}\line(1,0){60}}
\put(-20,30){\color{lightgray}\line(1,0){60}}
\put(10,-10){\color{lightgray}\line(-1,1){10}}
\put(0,0){\color{lightgray}\line(-1,5){10}}

\put(-40,-30){\color{lightgray}\line(1,0){60}}
\put(-30,-40){\color{cyan}\line(1,0){60}}
\put(-20,-50){\color{lightgray}\line(1,0){60}}
{\linethickness{1pt}
\put(-10,-30){\color{red}\line(1,-1){10}}
\put(10,-10){\color{red}\line(-1,-3){10}}
\put(-20,-10){\color{orange}\line(1,0){60}}
}

\put(0,0){\color{blue}\circle*{4}}
\put(0,0){\color{cyan}\circle*{2}}

\put(20,10){\circle*{3}}
\put(30,0){\color{blue}\circle*{4}}
\put(30,0){\color{cyan}\circle*{3}}
\put(30,0){\color{cyan}\circle*{3}}
\put(40,-10){\color{orange}\circle*{4}}
\put(43,-13){\color{orange}$L_1$}
\put(-10,10){\circle*{3}}
\put(-20,-10){\color{orange}\circle*{4}}
\put(-30,0){\color{blue}\circle*{4}}
\put(-30,0){\color{cyan}\circle*{3}}
\put(-40,10){\color{violet}\circle*{3}}
\put(-10,10){\circle*{3}}
\put(10,-10){\color{red}\circle*{4}}
\put(10,-10){\color{orange}\circle*{2}}
\put(11,-19){\color{red}$\vecbold{u}$}
\put(0,40){\color{cyan}\circle*{4}}
\put(1,42){\color{cyan}$\vecbold{v}$}
\put(20,50){\circle*{3}}
\put(30,40){\color{cyan}\circle*{3}}
\put(40,30){\circle*{3}}
\put(-10,50){\circle*{3}}
\put(-20,30){\circle*{3}}
\put(-30,40){\color{cyan}\circle*{3}}
\put(-40,50){\color{violet}\circle*{3}}
\put(-10,50){\circle*{3}}
\put(10,30){\circle*{3}}

\put(0,-40){\color{red}\circle*{4}}
\put(0,-40){\color{cyan}\circle*{2}}
\put(20,-30){\circle*{3}}
\put(30,-40){\color{cyan}\circle*{3}}
\put(40,-50){\circle*{3}}
\put(-20,-50){\circle*{3}}
\put(-30,-40){\color{cyan}\circle*{3}}
\put(-40,-30){\color{violet}\circle*{3}}
\put(-10,-30){\color{red}\circle*{4}}
\put(-11,-26){\color{red}$L_2$}
\put(10,-50){\circle*{3}}
\put(-20,1.5){\color{blue}$L$}
\end{picture}

Observe that $L_1\defeq\{(-1,1,0),(0,1,0),(1,1,0)\}$ and $L_2\defeq\{(0,-1,-1),(0,0,-1),(0,1,0)\}$ are two distinct lines in the plane $X=\overline{L\cup\{\vecbold{u}\}}$ that contain the point $\vecbold{u}$ and are disjoint with the line $L$. Also the lines $L_2$ and $L_3\defeq\{(-1,-1,-1),(-1,-1,0),(-1,-1,1)\}$ differ by the following properties: $\forall x\in L_3 \;(L_3\subseteq \overline{L_2\cup\{x\}})\quad\mbox{but}\quad \forall x\in L_2\;(L_2\not\subseteq \overline{L_3\cup\{x\}})$.

\end{proof}

\section{Ranked liners}

\begin{definition}\label{d:ranked} A liner $X$ is  called \defterm{ranked} if any two flats $A\subseteq B\subseteq X$ of the same finite rank $\|A\|=\|B\|<\w$ are equal.
\end{definition}

\begin{remark}\em A liner $X$ is ranked if and only if its rank function $\|\cdot\|$ is strictly monotone in the sense that distinct flats $A\subset B$ of finite rank have distinct ranks $\|A\|<\|B\|$.
\end{remark}

Definition~\ref{d:ranked} can be quantified as follows.

\begin{definition}\label{d:k-ranked} A liner $X$ is called \defterm{$\kappa$-ranked} for a cardinal $\kappa$ if any two flats $A\subseteq B\subseteq X$ of finite rank $\|A\|=\|B\|\le \kappa$ are equal.
\end{definition}

It is clear that a liner is ranked if and only if it is $\kappa$-ranked for every cardinal $\kappa$ if and only if it is $n$-ranked for every $n\in\w$. Every liner is $2$-ranked.

\begin{theorem}[{\cite[2.3.5]{Banakh}}]\label{t:ranked<=>EP} Let $\kappa$ be a cardinal. A liner is $\kappa$-ranked if and only if it has the $\kappa$-Exchange Property.
\end{theorem}

\begin{corollary}\label{c:ranked<=>EP} A liner is ranked if and only if it has the Exchange Property.
\end{corollary}

The notion of a ranked liner allows us to prove the following characterization of planes in $3$-ranked liners.

\begin{proposition}[{\cite[2.3.7]{Banakh}}]\label{p:SEP-plane}
Let $X$ be a $3$-ranked liner.
\begin{enumerate}
\item[\textup{(1)}] A set $A\subseteq X$ is a plane if and only if $A=\overline{\{x,y,z\}}$ for some points $x\in X$, $y\in X\setminus\{x\}$, and $z\in X\setminus\Aline xy$.
\item[\textup{(2)}] If a set $A\subseteq X$ is a plane, then $A=\overline{\{x,y,z\}}$ for every points $x\in A$, $y\in A\setminus\{x\}$, and $z\in X\setminus\Aline xy$.
\end{enumerate}
\end{proposition}

\begin{example} The liner in Example~\ref{ex:Tao} is not $3$-ranked: it is a plane containing a proper subplane.
\end{example}

\section{The relative rank and codimension of sets in liners}

\begin{definition} For two sets $A,B$ in a liner $X$, let $\|A\|_B$ be the smallest cardinality of a set $C\subseteq X$ such that $A\subseteq\overline{B\cup C}$. The cardinal $\|A\|_B$ is called the \defterm{$B$-relative rank} (or just the \defterm{$B$-rank}) of the set $A$ in $X$. If $B\subseteq A$, then the $B$-rank $\|A\|_B$ is denoted by $\dim_B(A)$ and called the \defterm{codimension} of $B$ in $A$.
\end{definition}

Observe that every set $A$ in a liner has rank $\|A\|=\|A\|_\emptyset$. So, the rank is a special case of the $B$-rank for $B=\emptyset$. 

It is easy to see that the relative rank is bimonotone in the sense that $\|A\|_C\le\|B\|_D$ for any sets $A\subseteq B$ and $D\subseteq C$ in a liner.

\begin{proposition}[{\cite[2.4.3]{Banakh}}]\label{p:rank-EP3} If a liner $X$ has the $\kappa$-Exchange Property for some cardinal $\kappa$, then for every subsets $A,B\subseteq X$ with $\|A\cup B\|\le \kappa$, there exist sets $A'\subseteq A$ and $B'\subseteq B\setminus\overline A$ such that $A\subseteq \overline {A'}$,
 $B\subseteq \overline{A'\cup B'}$, $|A'\cup B'|=\|A\cup B\|$, and $|B'|=\|B\|_A$.
\end{proposition}

\section{Independence in liners}

\begin{definition} A subset $A$ of a liner $X$ is \defterm{independent} if $a\notin\overline{A\setminus\{a\}}$ for every $a\in A$. 
\end{definition}

It is easy to see that any subset of an independent set in a liner is independent. 
The notion of an independent set is a partial case of a more general notion of a $B$-independent set.

\begin{definition} Let $B$ be a set in a liner $X$. A set $A\subseteq X$ is called \defterm{$B$-independent} if $a\notin \overline{B\cup(A\setminus\{a\})}$ for every $a\in A$.
\end{definition}

Therefore, a set $A$ in a liner $X$ is independent if and only if it is $\emptyset$-independent.

\begin{proposition}[{\cite[2.5.8]{Banakh}}]\label{p:union-of-independent} Let $A$ be a set in a liner $X$, $I$ be an $A$-independent set in $X$ and $J$ is an $(A\cup I)$-independent set in $X$. If the liner $X$ has the $\|A\cup I\cup J\|$-Exchange Property, then the set $I\cup J$ is $A$-independent.
\end{proposition}

\begin{proposition}[{\cite[2.5.9]{Banakh}}]\label{p:add-point-to-independent} Let $A,I$ be two sets in a liner $X$. If the set $I$ is $A$-independent and the liner $X$ is $\|A\cup I\|$-ranked, then for every point $x\in X\setminus\overline{A\cup I}$, the set $I\cup\{x\}$ is $A$-independent.
\end{proposition}

\begin{definition} Let $B,A$ be two sets in a liner $X$. A $B$-independent set $I\subseteq A$ is called a \defterm{maximal $B$-independent in} $A$ if $I=J$ for any $B$-independent set $J$ with $I\subseteq J\subseteq A$.
\end{definition} 

The Kuratowski--Zorn Lemma implies that for every sets $A,B$ in a liner $X$, every $B$-independent set $I\subseteq A$ can be enlarged to a maximal $B$-independent set $J\subseteq A$.

\begin{theorem}[{\cite[2.5.14]{Banakh}}]\label{t:Max=codim} Let $A,B$ be any sets in a liner $X$. If the liner $X$ has the $\|A\cup B\|$-Exchange Property, then the $A$-rank $\|B\|_A$ of the set $B$ is equal to the cardinality of any maximal $A$-independent set in $B$. 
\end{theorem}

\begin{corollary}[{\cite[2.5.15]{Banakh}}]\label{c:Max=dim} Let $A$ be a set in a liner $X$. If the liner $X$ has the $\|A\|$-Exchange Property, then the rank $\|A\|$ of the set $A$ is equal to the cardinality of any maximal indepen\-dent set in $A$. 
\end{corollary}

\begin{proposition}[{\cite[2.5.16]{Banakh}}]\label{p:rank+} For every sets $A,B,C$ in a liner $X$, 
$$\|C\|_A\le\|B\|_{A}+\|C\|_{B}.$$If $\overline A\subseteq \overline B\subseteq \overline C$ and the liner $X$ has the $\|C\|$-Exchange Property, then $\|C\|_A=\|B\|_A+\|C\|_B$.
\end{proposition}

\begin{corollary}[{\cite[2.5.17]{Banakh}}]\label{c:rank+} For every sets $A,B$ in a liner $X$, 
$\|B\|\le\|A\|+\|B\|_{A}$. If $\overline A\subseteq \overline B$ and the liner $X$ has the $\|B\|$-Exchange Property, then $\|B\|=\|A\|+\|B\|_A$.
\end{corollary}

\section{Regularity Axioms}

In this section we introduce three properties of liners, desribing the structure of the hull $\overline{A\cup\{a\}}$ of a flat $A$ with an attached singleton $\{a\}$.

For two subsets $A,B$ of a liner $X$ and a point $x\in X$, let 
$$\Aline Ax=\Aline xA\defeq\bigcup_{a\in A}\Aline ax\quad\mbox{and}\quad\Aline AB\defeq\bigcup_{x\in A}\bigcup_{y\in B}\Aline xy=\bigcup_{a\in A}\Aline aB=\bigcup_{b\in B}\Aline Ab.$$

\begin{definition}\label{d:regular} A liner $X$ is called
\begin{itemize}
\item \defterm{strongly regular} if for every  nonempty flat $A\subseteq X$ and point $b\in X\setminus A$, we have\newline $\overline{A\cup\{b\}}=\Aline Ab$;
\item \defterm{regular} if for every flat $A\subseteq X$ and points $a\in A $, $b\in X\setminus A$ we have\newline $\overline{A\cup\{b\}}=\bigcup_{y\in\Aline ab}\Aline Ay$;
\item \defterm{proregular} if for every flat $A\subseteq X$ and points $a\in A $, $b\in X\setminus A$ with $\Aline ab\ne\{a,b\}$, we have $\overline{A\cup\{b\}}=\bigcup_{y\in\Aline ab}\Aline Ay$;
\item \defterm{weakly regular} if for every flat $A\subseteq X$ and points $a\in A $, $b\in X\setminus A$ we have\newline $\overline{A\cup\{b\}}=\bigcup_{x\in A}\overline{\{a,b,x\}}$.
\end{itemize}
\end{definition}

Each strongly regular liner is regular and each regular liner is proregular and weakly regular. The notion of a (pro)regular liner can be quantified as follows.

\begin{definition}\label{d:k-proregular} Let $\kappa$ be a cardinal. A liner $X$ is called \defterm{$\kappa$-regular} (resp. \defterm{$\kappa$-proregular}\/) if for every set $A\subseteq X$ of cardinality $|A|<\kappa$ and every points $o\in\overline A $, $p\in X\setminus\overline A$ (such that $\Aline op\ne\{o,p\}$), we have $\overline{\{p\}\cup A}=\bigcup_{u\in\overline{o\,p}}\bigcup_{a\in \overline A}\Aline ua$.
\end{definition}

A liner is (pro)regular if and only if it is $\kappa$-(pro)regular for every cardinal $\kappa$. 

\begin{proposition}[{\cite[3.1.6]{Banakh}}]\label{p:k-regular<=>2ex} Let $\kappa$ be a cardinal. Every $\kappa$-proregular liner is $\kappa$-ranked and has the $\kappa$-Exchange Property.
\end{proposition}

\begin{example}\label{ex:Steiner13}The $13$-element group $\IZ_{13}$ endowed with the line relation

\centerline{$\Af\defeq\big\{(x,y,z)\in \IZ_{13}^3:\{x,y,z\}\in\big\{a+\{0,3,4\},a+\{0,5,7\}:a\in\IZ_{13}\big\}\big\}$} 

\noindent is an example of a finite liner, which is ranked but not $3$-proregular.
\end{example}

Proposition~\ref{p:k-regular<=>2ex} implies the following corollary.

\begin{corollary}[{\cite[3.1.8]{Banakh}}]\label{c:proregular=>ranked} Every proregular liner is ranked and has the Exchange Property.
\end{corollary}

The Exchange Property in weakly regular liners is equivalent to the $3$-Exchange Property.

\begin{proposition}[{\cite[3.1.9]{Banakh}}]\label{p:wr-ex<=>3ex} A weakly regular liner $X$ is ranked if and only if it $3$-ranked.
\end{proposition}

\begin{theorem}[{\cite[3.1.10+3.1.11]{Banakh}}]\label{p:reg<=>wreg+3-reg} A liner $X$ is regular if and only if it is $4$-regular if and only if it is weakly regular and $3$-regular.
\end{theorem}

\section{Modular and weakly modular liners}

\begin{theorem}[{\cite[5.1.1]{Banakh}}]\label{t:submodular} For any flats $A,B$ in an ranked liner $X$, we have the inequality
$\|A\cap B\|+\|A\cup B\|\le\|A\|+\|B\|$. 
If $\|A\cup B\|$ is infinite, then $\|A\cap B\|+\|A\cup B\|=\|A\|+\|B\|$.
\end{theorem} 

\begin{definition} A liner $X$ is called \textup{(}\defterm{weakly}\textup{)} \defterm{modular} if $\|A\cap B\|+\|A\cup B\|=\|A\|+\|B\|$ for any flats $A,B\subseteq X$ (with $A\cap B\ne\varnothing$).
\end{definition}

\begin{theorem}[{\cite[5.3.2]{Banakh}}]\label{t:modular} For a ranked liner $X$ the following conditions are equivalent:
\begin{enumerate}
\item[\textup{(1)}] $X$ is (weakly) modular;
\item[\textup{(2)}] $A\cap\overline{B\cup C}=\overline{(A\cap B)\cup C}$ for any flats $A,B\subseteq X$ and set $C\subseteq A$ (with $A\cap B\ne\varnothing$);
\item[\textup{(3)}] $A\cap\overline{B\cup C}=\overline{(A\cap B)\cup C}$ for any flats $A,B\subseteq X$ and set $C\subseteq A$ with $\|A\cup B\|<\w$ (and $A\cap B\ne\varnothing$).
\end{enumerate}
\end{theorem}

Two lines $A,B$ in a liner are \defterm{coplanar} if they are contained in some plane.

\begin{theorem}[{\cite[5.5.1]{Banakh}}]\label{t:modular<=>}
A liner $X$ is modular if and only if it is strongly regular if and only if $X$ contains no disjoint coplanar lines.
\end{theorem}

\begin{theorem}[{\cite[5.4.1]{Banakh}}]\label{t:w-modular<=>}For a liner $X$ the following conditions are equivalent:
\begin{enumerate}
\item[\textup{(1)}] $X$ is weakly modular;
\item[\textup{(2)}] $X$ is ranked and weakly regular;
\item[\textup{(3)}] $X$ is ranked and for every planes $P\cup\Pi$ with $\|P\cup\Pi\|=4$, the intersection $P\cap\Pi$ is not a singleton.
\end{enumerate}
\end{theorem}

\begin{corollary} Every regular liner is weakly modular.
\end{corollary}

\begin{example} The Euclidean plane is an example of a regular weakly modular liner which is not strongly regular and hence not modular.
\end{example}

\section{Parallelity Postulates and Axioms}

In this section we introduce many parallelity postulates and parallelity axioms that will be studied in the next sections. We start with four classical Parallelity Postulates.

\begin{definition}\label{d:PPBL} A liner $X$ is defined to be
\begin{itemize}
\item \defterm{Proclus} 
if for every plane $P\subseteq X$, line $L\subseteq P$ and point $x\in P\setminus L$ there exists at most one line $\Lambda$ in $X$ such that $x\in \Lambda\subseteq P\setminus L$;
\item \defterm{Playfair} 
if for every plane $P\subseteq X$, line $L\subseteq P$ and point $x\in P\setminus L$ there exists a unique line $\Lambda$ in $X$ such that $x\in \Lambda\subseteq P\setminus L$;
\item \defterm{Bolyai} 
if for every plane $P\subseteq X$, line $L\subseteq P$ and point $x\in P\setminus L$ there exists at least one line $\Lambda$ in $X$ such that $x\in \Lambda\subseteq P\setminus L$;
\item \defterm{Lobachevsky} 
if for every plane $P\subseteq X$, line $L\subseteq P$ and point $x\in P\setminus L$ there exist at least two distinct lines $\Lambda$ in $X$ such that $x\in \Lambda\subseteq P\setminus L$.
\end{itemize}
\end{definition}

These four definitions are partial cases of the following Parallel Postulates involving a cardinal parameter $\kappa$.

\begin{definition}\label{d:k-PPB} Let $\kappa$ be a cardinal. A liner $X$ is defined to be 
\begin{itemize}
\item \defterm{$\kappa$-hypoparallel} if for every plane $P\subseteq X$, line $L\subseteq P$ and point $x\in P\setminus L$ there exist at most $\kappa$ lines $\Lambda$ such that $x\in \Lambda\subseteq P\setminus L$;
\item \defterm{$\kappa$-parallel} if for every plane $P\subseteq X$, line $L\subseteq P$ and point $x\in P\setminus L$ there exist exactly $\kappa$ lines $\Lambda$ such that $x\in \Lambda\subseteq P\setminus L$;
\item \defterm{$\kappa$-hyperparallel} if for every plane $P\subseteq X$, line $L\subseteq X$ and point $x\in P\setminus L$ there exist at least $\kappa$  lines $\Lambda$ such that $x\in \Lambda\subseteq P\setminus L$.
\end{itemize}
\end{definition}

It is clear that a liner is $\kappa$-parallel if and only if it is both $\kappa$-hypoparallel and $\kappa$-hyperparallel.

Observe that a line $X$ is Proclus, Playfair, Bolyai, Lobachevsky if and only if it is $1$-hypoparallel, $1$-parallel, $1$-hyperparallel, $2$-hyperparallel, respectively.

 %
%

In Theorems~\ref{t:Proclus<=>} and \ref{t:Playfair<=>}, we shall show that  Proclus and Playfair Parallel Postulates are first-order properties of liners, whose first-order characterizations involve the following parallelity axioms.

\begin{definition}\label{d:PP} A liner $X$ is defined to be
\begin{itemize}
\item \defterm{projective} if $\forall o,x,y\in X\;\;\forall p\in\Aline xy\;\;\forall v\in\Aline oy\setminus\{p\}\;\;(\Aline vp\cap\Aline ox\ne \varnothing)$;
\item \defterm{proaffine} if $\forall o,x,y\in X\;\forall p\in\Aline xy\setminus\Aline ox\;\exists u\in\Aline oy\;\forall v\in\Aline oy\setminus \{u\}\;\;(\Aline vp\cap\Aline ox\ne \varnothing)$;
\item \defterm{affine} if $\forall o,x,y\in X\;\forall p\in\Aline xy\setminus\Aline ox\;\exists u\in\Aline oy\;\forall v\in\Aline oy\;\;(u=v\;\Leftrightarrow\;\Aline vp\cap\Aline ox=\overline \varnothing)$;
\item \defterm{hyperaffine} if $\forall o,x,y\in X\;\forall p\in\Aline xy\setminus\Aline ox\;\exists u\in\Aline oy\;\;\;(\Aline up\cap\Aline ox= \varnothing)$;
\item \defterm{hyperbolic} if $\forall o,x,y\in X\,\forall p\in\Aline xy\setminus(\Aline ox\,{\cup}\,\Aline oy)\,\exists u,\!v{\in}\Aline oy\;(u\ne v\,\wedge\,\Aline up\cap\Aline ox{=}\varnothing{=}\Aline up\cap\Aline ox)$;
\item \defterm{injective} if $\forall o,x,y\in X\;\forall p\in\Aline xy\setminus(\Aline ox\,{\cup}\,\Aline oy)\;\forall v\in\Aline oy\setminus\{o,y\}\;\;(\Aline vp\cap\Aline ox=\varnothing)$.
\end{itemize} The first-order formulas defining projective, proaffine, affine, hyperaffine, and hyperbolic liners will be called the \defterm{projectivity}, \defterm{proaffinity}, \defterm{affinity}, \defterm{hyperaffinity}, \defterm{hyperbolicity}, and \defterm{injectivity axioms}, respectively. 
\end{definition}

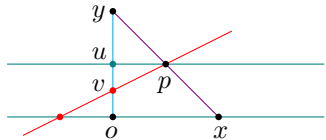
\begin{picture}(100,55)(-150,-10)

\put(0,0){\color{teal}\line(1,0){120}}
\put(40,0){\color{cyan}\line(0,1){40}}
\put(40,40){\color{violet}\line(1,-1){40}}
\put(0,20){\color{teal}\line(1,0){120}}
\put(60,20){\color{red}\line(2,1){25}}
\put(60,20){\color{red}\line(-2,-1){54}}

\put(40,0){\circle*{2.5}}
\put(37,-8){$o$}
\put(80,0){\circle*{2.5}}
\put(78,-8){$x$}
\put(40,40){\circle*{2.5}}
\put(32,38){$y$}
\put(40,20){\color{teal}\circle*{2.5}}
\put(32,22){$u$}
\put(40,10){\color{red}\circle*{2.5}}
\put(32,10){$v$}
\put(60,20){\circle*{2.5}}
\put(57,11){$p$}
\put(20,0){\color{red}\circle*{2.5}}

\end{picture}

\begin{remark} The proaffinity, affinity, hyperaffinity, and hyperbolicity axioms are first order counterparts of the parallel postulates of Proclus, Playfair, Bolyai, and Lobachevsky, respectively. The projectivity axiom has been explicitely formulated  Oswald Veblen 
in his axiom system for Projective Geometry. The injectivity axiom (as an antonym to the projectivity) is well-known in the theory of unitals, as the absence of O'Nan configurations (known also as Pasch configurations), see Section~16 for more information on the classical unitals.
\end{remark}

\begin{remark} Every $3$-long  injective liner is hyperaffine, and every $4$-long injective liner is hyperbolic.
\end{remark}

The above parallelity postulates and axioms relate as follows.
$$
\xymatrix@C=20pt@R=20pt{
\mbox{$0$-hypoparallel}\ar@{=>}[r]\ar@{<=>}[d]&\mbox{Proclus}\ar@{<=>}[d]&\mbox{Playfair}\ar@{=>}[l]\ar@{<=>}[d]\ar@{=>}[r]&
\mbox{Bolyai}\ar@{<=>}[d]&\mbox{Lobachevsky}\ar@{=>}[l]\ar@{<=>}[d]\\
\mbox{$0$-parallel}\ar@{=>}[r]\ar@{<=>}[d]&\mbox{$1$-hypoparallel}\ar@{=>}[d]&\mbox{$1$-parallel}\ar@{=>}[l]\ar@{=>}[d]\ar@{=>}[r]&
\mbox{$1$-hyperparallel}&\mbox{$2$-hyperparallel}\ar@{=>}[l]\\
\mbox{projective}\ar@{=>}[r]&\mbox{proaffine}&\mbox{affine}\ar@{=>}[l]\ar@{=>}[r]&\mbox{hyperaffine}\ar@{=>}[u]&\mbox{hyperbolic}\ar@{=>}[l]\ar@{=>}[u]
}
$$

\begin{proposition}[{\cite[3.2.8]{Banakh}}] A liner $X$ is injective if and only if among any four distinct lines in $X$ two are disjoint.
\end{proposition}

\begin{example} For every $n\ge 2$, the subliner $B=\{x\in \IR^n:\sum_{i\in n}x_i^2<1\}$ of the Euclidean space $\IR^n$ is hyperbolic. This hyperbolic liner is called the \defterm{Beltrami--Klein model} of hyperbolic geometry.
\end{example}

\section{Balanced liners}

\begin{definition} Given two cardinals $\kappa,\mu$, we say that a liner $X$ is \defterm{$\kappa$-balanced} with $|X|_\kappa=\mu$ if every flat $A\subseteq X$ of rank $\|A\|=\kappa$ has cardinality $|A|=\mu$. A liner $X$ is \defterm{balanced} if it is $\kappa$-balanced for every cardinal $\kappa$.
\end{definition}

 For a $\kappa$-balanced liner $X$, the cardinal $|X|_\kappa$ is well-defined if $X$ contains a flat $F\subseteq X$ of rank $\|F\|=\kappa$. In this case $|X|_\kappa=|F|$ for every flat $F\subseteq X$ of rank $\|F\|=\kappa$.

It is clear that every liner $X$ is $1$-balanced with $|X|_1=1$. A liner $X$ is $2$-balanced if and only if all lines have the same cardinality $|X|_2$.  

\begin{theorem}[{\cite[3.3.3]{Banakh}}]\label{t:2-balanced-dependence} Let $X$ be a $2$-balanced liner and $\mathcal L$ be the family of all lines in $X$. If $|X|_2<|X|$, then for every point $x\in X$, the cardinality of the family $\mathcal L_x\defeq\{L\in\mathcal L:x\in L\}$ satisfies two equations
$|X|-1=(|X|_2-1)\cdot |\mathcal L_x|$ and $|X|_2\cdot|\mathcal L|=|X|\cdot|\mathcal L_x|$, and hence $|\mathcal L_x|$ does not depend on $x$. If $X$ is finite, then $|X|_2$ divides the number $(|\mathcal L_x|-1)\cdot|\mathcal L_x|$.
\end{theorem}

\begin{remark} Theorem~\ref{t:2-balanced-dependence} implies that the cardinality $v\defeq|X|$ of any finite $2$-balanced liner $X$ with $k\defeq |X|_2<|X|$ satisfies the equality $v=1+(k-1)\cdot r$ for some number $r\ge k$ such that $k$ divides $(r-1){\cdot}r$ (denoted by $k\,|\,(r-1)r)$. The latter divisibility condition implies (and in fact, is equivalent to) the divisibilities $(k-1)\,|\,(v-1)$ and $(k-1)k\,|\,(v-1)v$. Pairs of positive integer numbers $(k,v)$ satisfying the divisibility conditions 
\smallskip

\centerline{$(k-1)\,|\,(v-1)\quad\mbox{and}\quad(k-1)k\,|\,(v-1)v$}
\smallskip

\noindent are called \defterm{admissible}. Therefore, for every finite $2$-balanced liner $X$, the pair $(|X|_2,|X|)$ is admissible. The admissibility of a pair $(k,v)$ does not imply the existence of a $2$-balanced liner $X$ with $(|X|_2,|X|)=(k,v)$. 
For example, the pairs $(6,36)$ and $(7,43)$ are admissible, but liners $X$ with  $(|X|_2,|X|)\in\{(6,36),(7,43)\}$ do not exist, by the Bruck--Ryser's Theorem~\ref{t:Bruck-Ryser}. The pair $(6,46)$ is admissible but no liner $X$ with $(|X|_2,|X|)=(6,46)$ exists, by a result of Houghten, Thiel, Janssen,  and Lam \cite{HTJL2001}, proved by a computer search.  On the other hand, the admissibility condition is asymptotically sufficient for the existence of finite $2$-balanced liners, according to the following fundamental theorem, proved by Wilson in \cite{Wilson1}, \cite{Wilson2}, \cite{Wilson3}.
\end{remark}

\begin{theorem}[Wilson, 1975]\label{t:Wilson} For every integer numbers $k$ there exists a number $v_k$ such that for every admissible pair $(k,v)$ with $v>v_k$, there exists a $2$-balanced liner $X$ such that $|X|_2=k$ and $|X|=v$. 
\end{theorem}

\begin{remark} We can assume that the number $v_k$ in Wilson's Theorem~\ref{t:Wilson} is the smallest possible. In this case it is determined uniquely. The Wilson proof of Theorem~\ref{t:Wilson} does not provide an explicit formula for finding the number $v_k$, it just claims that it does exist. However, for $k\le 9$ the following upper and lower bounds for the numbers $v_k$ are known (see \cite[\S II.3.1]{HCD}):
$$
\begin{array}{|c|cccccccc|}
\hline
k&2&3&4&5&6&7&8&9\\
\hline
v_k\ge&2&3&4&5&46&43&8&9\\
v_k\le &2&3&4&5&801&2605&3753&16497\\
\hline
\end{array}
$$As we have already mentioned, liners $X$ with $|X|_2=6$ and $X\in\{36,46\}$ do not exist. The existence of liners $X$ with $|X|_2=6$ is undecided for the following values $v$ (see Table 3.4 in \cite[II.3.1]{HCD} and \cite{BHR2025}): 51, 61, 81, 166, 231, 256, 261, 286, 316, 321, 346, 351, 376, 406, 411, 436, 471, 501, 561, 591,
616, 646, 651, 676, 771, 796, 801.
\end{remark}
 
\begin{theorem}[{\cite[3.3.7]{Banakh}}]\label{t:2-balance+k-parallel=>3-balance} If a liner $X$ is $2$-balanced and $\kappa$-parallel for some cardinal $\kappa$, then $X$ is $3$-balanced with $|X|_3=1+(\kappa+|X|_2)(|X|_2-1)$. If the cardinal $|X|_3$ is well-defined and finite, then the liner $X$ is $3$-ranked.
\end{theorem}

\begin{proposition}[{\cite[3.3.8]{Banakh}}]\label{p:23-balance=>k-parallel} If a liner $X$ of rank $\|X\|\ge 3$ is $2$-balanced and $3$-balanced, then there exists a unique cardinal $\lambda\ge|X|_2$ such that $|X|_3-1=\lambda\cdot(|X|_2-1)$. If $|X|_2<|X|_3$, then $X$ is $\kappa$-parallel for the unique cardinal $\kappa$ such that $\lambda=\kappa+|X|_2$.
\end{proposition} 

\begin{theorem}[{\cite[3.3.9]{Banakh}}]\label{t:k-parallel=>2-balance} Assume that a liner $X$ is $\kappa$-parallel for some finite cardinal $\kappa$. If $\kappa>0$ or $X$ is $3$-long, then $X$ is $2$-balanced, $3$-balanced and $|X|_3=1+(|X|_2-1)(\kappa+|X|_2)$. If $|X|_2$ is finite, then $|X|_2$ divides $(\kappa-1)\kappa$.
\end{theorem}

Theorems~\ref{t:2-balance+k-parallel=>3-balance} and Proposition~\ref{p:23-balance=>k-parallel} imply the following characterization.


\begin{corollary}[{\cite[3.3.11]{Banakh}}] For a $3$-long finite liner $X$ of rank $\|X\|\ge 3$, the following conditions are equivalent:
\begin{enumerate}
\item[\textup{(1)}] $X$ is $\kappa$-parallel for some cardinal $\kappa$;
\item[\textup{(2)}] $X$ is $2$-balanced and $3$-balanced;
\item[\textup{(3)}] $X$ is $2$-balanced, $3$-balanced and $\kappa$-parallel for a unique cardinal $\kappa$\newline  satisfying the equation $|X|_3-1=(\kappa+|X|_2)(|X|_2-1)$.
\end{enumerate}
\end{corollary}


\begin{proposition}[{\cite[3.3.13]{Banakh}}]\label{p:2-balanced=>w-balanced} If a liner $X$ is $2$-balanced, then for every infinite cardinal $\kappa$, the liner $X$ is $\kappa$-balanced with $|X|_\kappa=\max\{\kappa,|X|_2\}$.
\end{proposition}

\begin{theorem}[{\cite[3.3.14]{Banakh}}]\label{t:wr+k-parallel=>n-balanced} If a weakly regular $3$-ranked liner $X$ is $2$-balanced and $\kappa$-parallel for some cardinal $\kappa$, then $X$ is balanced with $|X|_{n}=1+(|X|_2-1)\sum_{r=0}^{n-2}(\kappa+|X|_2-1)^r$ for every finite cardinal $n\ge 2$.
\end{theorem}

\begin{question} For which cardinals $\kappa$ there exist weakly regular $3$-long $\kappa$-balanced ranked liners of rank $\|X\|>3$?
\end{question} 

\begin{proposition}[{\cite[3.3.16]{Banakh}}] Every $2$-balanced liner $X$ with $|X|_2{<}|X|{<}\,|X|_2^3{-}2{\cdot}|X|_2^2{+}2{\cdot}|X|_2$ is a ranked plane. 
\end{proposition}

\begin{theorem}[{\cite[3.3.17]{Banakh}}]\label{t:affine=>Avogadro} Every affine liner $X$ is $2$-balanced.
\end{theorem}

\begin{remark}\label{r:order-affine} By Theorem~\ref{t:affine=>Avogadro}, all lines in an affine liner have the same cardinality, which is called the \defterm{order} of the affine liner. By the Veblen--Young Theorem (1908), the order of any finite affine regular liner of dimension $\ge 3$ is a power of a prime number. There is a (still unproved) conjecture that the same result is true for finite affine planes: their orders are powers of primes. 
\end{remark}


\begin{proposition}[{\cite[3.3.24]{Banakh}}]\label{p:Avogadro-projective} Every $3$-long projective liner $X$ is $2$-balanced.
\end{proposition}

If a liner $X$ of rank $\|X\|\ge 2$ is not $3$-long, then it contains a line of length $2$ and there is no further restrictions on the length of lines  in such liners.

\begin{example} For every set $K$ of non-zero cardinals with $1,2\in K$, there exists a projective liner $X$ of cardinality $|X|=\sum_{k\in K}k$ and rank $2\cdot|K|-1$ such that $K=\{|\Aline xy|:x,y\in X\}$.
\end{example}

\begin{corollary}[{\cite[3.3.26+8.5.2]{Banakh}}]\label{c:L<L+2} Any lines $L,\Lambda$ in a $3$-long proaffine (and Proclus) liner $X$ have $|L|\le|\Lambda|+2$ (and $|L|\le|\Lambda|+1$). 
\end{corollary}

\begin{corollary}[{\cite[3.3.22]{Banakh}}] Every $\w$-long proaffine liner is $2$-balanced.
\end{corollary}

\begin{remark} If all lines in a projective liner have the same finite cardinality $\lambda$, then the number $\lambda-1$ is called the \defterm{order} of the projective liner. By the Veblen--Young Theorem, the order of a finite projective liner of dimension $\ge 3$ is a power of a prime number. There is a (still unproved) conjecture that the same result holds for projective planes: their order is a power of a prime number. Using projective completions (discussed in Chapter~7 of \cite{Banakh}), it can be shown that a  projective plane of order $n$ exists if and only if an affine plane of order $n$ exists.
\end{remark}

The unique available result on possible orders of finite affine planes is the following fundamental theorem of Bruck--Ryser \cite{BR1949}.

\begin{theorem}[Bruck--Ryser, 1949; {\cite[25.5.1]{Banakh}}]\label{t:Bruck-Ryser} If a number $n\in (4\IZ+1)\cup(4\IZ+2)$ is an order of an affine or projective plane, then $n=x^2+y^2$ for some integer numbers $x,y$. 
\end{theorem}

Bruck--Ryser Theorem implies that the order of an affine or projective plane cannot be equal to $6$, $14$, or $22$. On the other hand, the number $10=4\cdot 2+2=1^2+9^2$ is the sum of two squares but no affine plane of order 10 exists, as was shown by Lam, Thiel and Swiercz~\cite{LTS1989} in 1989, using heavy computer calculations. The problem of the (non-)existence of an affine or projective plane of order $12$ remains open, out of reach of modern (super)computers.

\section{The regularity of balanced liners}

In this section we discuss the structure of weakly regular balanced liners.

\begin{lemma}[{\cite[5.6.1]{Banakh}}]\label{l:P=union-of-lines} Let $X$ be a weakly modular liner of rank $\|X\|\ge 4$. If some plane $P\subseteq X$ contains two disjoint lines, then $P$ is the union of a family of pairwise disjoint lines.
\end{lemma}

\begin{corollary}[{\cite[5.6.3]{Banakh}}]\label{c:x2|x3} Let $X$ be a balanced weakly modular liner of rank $\|X\|\ge 4$. If $X$ is not projective, then $|X|_3=|X|_2\cdot\kappa$ for some cardinal $\kappa\ge |X|_2$.
\end{corollary}

\begin{lemma}[{\cite[5.6.4]{Banakh}}]\label{l:disjoint-planes} Let $X$ be a $2$-balanced and $3$-balanced weakly modular liner. If $|X|_3<\w$ and $X$ is not projective, then for every flat $Y\subseteq X$ of rank $\|Y\|=4$, every plane $P\subseteq Y$ and every point $y\in Y\setminus P$, there exists a plane $\Pi\subseteq Y$ such that $y\in\Pi\subseteq Y\setminus P$.
\end{lemma}

The following lemma is a higher-dimensional counterpart of Lemma~\ref{l:P=union-of-lines}.

\begin{lemma}[{\cite[5.6.5]{Banakh}}]\label{l:Y=union-of-planes} Let $X$ be a weakly modular liner of rank $\|X\|\ge 5$. If some flat $Y\subseteq X$ of rank $\|Y\|=4$  contains two disjoint planes, then $Y$ is the union of a family of pairwise disjoint planes.
\end{lemma}

This lemma implies the following higher-dimensional counterpart of Corollary~\ref{c:x2|x3}.

\begin{corollary}[{\cite[5.6.7]{Banakh}}]\label{c:x3|x4} Let $X$ be a balanced weakly modular liner of rank $\|X\|\ge 5$. If $X$ is not projective, then $|X|_4=|X|_3\cdot\kappa$ for some cardinal $\kappa\ge |X|_2$.
\end{corollary}

The principal results of this section are the following theorems of Jean Doyen and Xavier Hubault  \cite{DH1971}.

\begin{theorem}[Doyen--Hubaut, 1971; {\cite[5.6.8+5.6.10]{Banakh}}]\label{t:DH1} Let $X$ be a balanced weakly modular liner of rank $\|X\|\ge 4$. If $|X|_3<\w$, then $X$ is $p$-parallel for some $p\in\{0,1,|X|_2^2\,{-}\,|X|_2\,{+}\,1,\break |X|^3_2+1\}$. If $\|X\|\ge 5$, then $p\in\{0,1\}$ and hence $X$ is proaffine.
\end{theorem}


\begin{problem} Is every finite balanced weakly regular liner of rank $\ge 4$ proaffine?
\end{problem}

\begin{theorem}[{\cite[5.6.12]{Banakh}}]\label{t:balanced<=>ranked} For a balanced liner $X$ with $\|X\|\ge 4$ and $|X|_3<\w$, the following conditions  are equivalent:
\begin{enumerate}
\item[\textup{(1)}]  $X$ is regular;
\item[\textup{(2)}]  $X$ is $p$-parallel for some $p\in\{0,1\}$;
\item[\textup{(3)}]  $X$ is projective or Playfair;
\item[\textup{(4)}]  $X$ is Proclus;
\item[\textup{(5)}]  $X$ is $3$-proregular and proaffine.
\end{enumerate}
If $\|X\|\ge5$, then the conditions \textup{(1)--(5)} are equivalent to:
\begin{enumerate}
\item[\textup{(6)}] $X$ is weakly regular;
\item[\textup{(7)}] $X$ is weakly modular.
\end{enumerate}
\end{theorem}

\begin{corollary}[{\cite[5.6.13]{Banakh}}] For a finite liner $X$ of rank $\|X\|\ge 5$, the following conditions are equivalent:
\begin{enumerate}
\item[\textup{(1)}] $X$ is regular and Playfair;
\item[\textup{(2)}] $X$ is balanced, weakly regular and not strongly regular.
\item[\textup{(3)}] $X$ is balanced, weakly modular and not modular.
\end{enumerate}
\end{corollary}

\begin{problem} Is every finite balanced finite affine liner $X$ of rank $4$ regular?
\end{problem}

\begin{remark} Computer calculations show every balanced regular liner $X$ with $|X|_2\le 4$ is proaffine.
\end{remark}

\begin{example}\label{ex:Denniston} There exists a hyperbolic  regular balanced liner $X$ with $|X|_2=8$ and $|X|=|X|_3=120$. This liner is one of maximal arcs, discovered by Denniston \cite{Denniston1969}.
\end{example}

\section{Projective liners}

In this section we present some characterizations of projective liners. 

Two lines $A,B$ in a liner $X$ are \defterm{skew} if $\|A\cup B\|=4$. It is clear that any skew lines are disjoint. The converse is true in projective liners.

\begin{theorem}[{\cite[3.4.1+5.5.1]{Banakh}}]\label{t:projective<=>} For every liner $X$, the following conditions are equivalent:
\begin{enumerate}
\item[\textup{(1)}] $X$ is $0$-hypoparallel;
\item[\textup{(2)}] $X$ is $0$-parallel;
\item[\textup{(3)}] any disjoint lines in $X$ are skew;
\item[\textup{(4)}] $X$ projective;
\item[\textup{(5)}] $X$ is strongly regular;
\item[\textup{(6)}] $X$ is modular.
\end{enumerate}
\end{theorem}



\begin{corollary}[{\cite[3.4.3]{Banakh}}]\label{c:Steiner-projective<=>} A $2$-balanced liner $X$ with $|X|_2<\w$ is projective if and only if $X$ is $3$-balanced and $|X|_3=|X|_2^2-|X|_2+1$.
\end{corollary}

Theorems~\ref{t:projective<=>} and \ref{t:wr+k-parallel=>n-balanced}  imply the following corollary.

\begin{corollary}[{\cite[3.4.4]{Banakh}}]\label{c:projective-order-n} Every finite projective liner $X$ of order $\ell$ contains exactly\newline $\sum_{k=0}^{\|X\|-1}\ell^k=\frac{\ell^{\|X\|}-1}{\ell-1}$ points.
\end{corollary}

By Proposition~\ref{p:Avogadro-projective}, every $3$-long projective liner is $2$-balanced.  Now we describe the structure of projective liners which are not $3$-long. This description involves the notion of a maximal $n$-long flat in a liner. 

\begin{definition}\label{d:max-long} Let $\kappa$ be a cardinal number. A flat $B$ in a liner $X$ is called \defterm{$\kappa$-long} if $|\Aline xy|\ge \kappa$ for every distinct points $x,y\in B$. A $\kappa$-long flat $B$ in $X$ is called \defterm{maximal $\kappa$-long} if every $\kappa$-long flat $A\subseteq X$ containing $B$ coincides with $B$.
\end{definition}

\begin{remark} By Definition~\ref{d:max-long}, every flat $F\subseteq X$ of cardinality $|F|\le 1$ in a liner $X$ is $\kappa$-long for every cardinal $\kappa$. The Kuratowski-Zorn Lemma implies that every $\kappa$-long flat in a liner $X$ is a subset of a maximal $\kappa$-long flat in $X$. In particular, every point $x$ of a liner $X$ belongs to some maximal $\kappa$-long flat (which can be equal to the singleton $\{x\}$, if the cardinal $\kappa$ is too large).
\end{remark}


A family of sets $\F$ is called \defterm{disjoint} if $A\cap B=\varnothing$ for every distinct sets $A,B\in\F$. The following theorem shows that the line structure of a projective liner is uniquely determined by the structures of its maximal $3$-long flats.

\begin{theorem}[{\cite[3.4.11]{Banakh}}]\label{t:projective=+3long} The family $\mathcal M$ of maximal $3$-long flats in a projective liner $X$ is disjoint and $\Aline ab=\{a,b\}$ for any distinct flats $A,B\in \mathcal M$ and points $a\in A$ and $b\in B$.
\end{theorem}

Next, we describe the structure of projective planes which are not $3$-long.

\begin{definition} A liner $X$ is called a \defterm{near-pencil} if $X$ contains a line $L\subseteq X$ such that $X\setminus L$ is a singleton.
\end{definition}

\noindent In the following theorem by a \defterm{projective plane} we understand a projective liner of rank $3$. 

\begin{theorem}[{\cite[3.4.13]{Banakh}}] A liner is a near-pencil if and only if it is projective plane which is not $3$-long.
\end{theorem}



\begin{corollary}[{\cite[3.4.16]{Banakh}}]\label{c:line-meets-hyperplane} Every line in a projective liner $X$ has nonempty intersection with every hyperplane in $X$.
\end{corollary}

We recall that a flat $H$ in a liner $X$ is a \defterm{hyperplane} if $H=X$ and every flat $F$ in $X$ with $H\subseteq F\subseteq X$ is equal to $H$ or $X$.

\begin{theorem}[Bruck, 1955; {\cite[3.5.1]{Banakh}}]\label{t:Bruck55} Let $X$ be a projective liner of finite order $n$ and $P$ be a projective subliner of order $p<n$ in $X$. If $\|P\|\ge 3$, then either $n=p^2$ or else $n\ge p^2+p$.
\end{theorem}

Theorem~\ref{t:Bruck55} motivates the following long-standing open problem.

\begin{problem} Is there a projective plane $X$ of composite order $p^2+p$ for some $p$? In particular, is there a projective plane of order $12=3^2+3$?
\end{problem}

\section{Proaffine and Proclus liners} 

In this section we study the interplay between proaffine and Proclus liners, and will prove that Proclus liners are exactly proaffine $3$-proreglar liners.

\begin{theorem}[{\cite[3.6.1]{Banakh}}]\label{t:Proclus<=>}For a liner $X$, the following conditions are equivalent:
\begin{enumerate}
\item[\textup{(1)}] $X$ is Proclus;
\item[\textup{(2)}] $X$ is proaffine and $3$-proregular;
\item[\textup{(3)}]  for every line $L\subseteq X$ and points $o\in L$, $x\in X\setminus L$,  $y\in\overline{L\cup\{x\}}$,\newline the set $\{v\in\Aline ox:\Aline vy\cap L=\varnothing\}$ contains at most one point.
\end{enumerate}
\end{theorem}

The following proposition is a version of the classical Proclus Parallelity Postulate.

\begin{proposition}[Proclus Postulate; {\cite[3.6.10]{Banakh}}]\label{p:Proclus-Postulate} Let $P$ be a plane in a Proclus liner $X$ and $L,L',\Lambda$ be three lines in the plane $P$. If $L\cap L'=\varnothing$ and $|\Lambda\cap L|=1$, then $|\Lambda\cap L'|=1$.
\end{proposition}

\begin{corollary}[{\cite[3.6.11]{Banakh}}]\label{c:Proclus-par=} Let $L,\Lambda$ be two disjoint lines in a Proclus liner $X$. If $\|L\cup\Lambda\|\le 3$, then $|L|=|\Lambda|$.
\end{corollary}

\section{Affine and Playfair liners}

In this section we study the interplay between affine and Playfair liners, and show that Playfair liners are exactly affine $3$-regular $3$-long liners.

\begin{theorem}[{\cite[3.7.1]{Banakh}}]\label{t:affine-char1} For a liner $X$, the following conditions are equivalent:
\begin{enumerate}
\item[\textup{(1)}] $X$ is affine and $3$-proregular;
\item[\textup{(2)}] $X$ is affine and $3$-regular;
\item[\textup{(3)}] for every line $L\subseteq X$ and points $o\in L$, $y\in X\setminus L$,  $z\in\overline{L\cup\{y\}}\setminus L$,\newline the set $\{u\in\Aline oy:\Aline uz\cap L=\emptyset\}$ is a singleton;
\item[\textup{(4)}] for every line $L\subseteq X$ and point $p\in X\setminus L$ with $\overline{L\cup\{p\}}\ne L\cup\{p\}$,\newline there exists a unique line $\Lambda$ such that $p\in \Lambda\subseteq\overline{\{p\}\cup L}\setminus L$.
\end{enumerate}
\end{theorem}

\begin{theorem}[{\cite[3.7.2]{Banakh}}]\label{t:Playfair<=>} For a liner $X$ of cardinality $|X|>2$, the following conditions are equivalent:
\begin{enumerate}
\item[\textup{(1)}] $X$ is Playfair;
\item[\textup{(2)}]  $X$ is affine, $3$-regular, and $3$-long;
\item[\textup{(3)}]  for every line $L\subseteq X$ and point $x\in X\setminus L$,\newline there exists a unique line $\Lambda$ such that $x\in \Lambda\subseteq\overline{L\cup\{x\}}\setminus L$.
\end{enumerate}
\end{theorem}

\begin{corollary}[{\cite[3.7.3]{Banakh}}]\label{c:affine-cardinality} Every $3$-long affine regular liner $X$ of finite rank $r$ is $2$-balarnced and has cardinality $|X|=|X|_2^{r-1}$.
\end{corollary}  

\begin{corollary}[{\cite[3.7.4]{Banakh}}]\label{c:Playfair<=>23balanced} A liner $X$ with finite lines is Playfair if and only if $X$ is $2$-balanced, $3$-balanced and $|X|_3=|X|_2^2$.
\end{corollary}

\begin{proposition}[{\cite[3.7.5]{Banakh}}]\label{p:Playfair-plane<=>} A liner $X$ is a Playfair plane if and only if $X$ is $3$-long, $\|X\|>2$, and for every line $L\subseteq X$ and point $x\in X\setminus L$ there exists a unique line $L_x$ in $X$ such that $x\in L_x\subseteq X\setminus L$.
\end{proposition}




\begin{theorem}[{\cite[4.1.1]{Banakh}}]\label{t:4-long-affine} Every $4$-long affine liner is Playfair and regular.
\end{theorem}

The following examples show that Theorem~\ref{t:4-long-affine} does not extend to $3$-long affine liners and also to $4$-long Proclus liners.

\begin{example}[{\cite[3.7.7+4.5.7]{Banakh}}] The ring $\IZ_{15}$ endowed with the family of lines
\smallskip

\centerline{$\mathcal L=\{x+L:x\in\IZ_{15},\;L\in\mathcal B\}\mbox{ \  where \ }\mathcal B=\big\{\{0,1,4\},\{0,6,8\},\{0,5,10\}\big\}$}
\smallskip

\noindent is an example of a $5$-parallel non-regular ranked balanced affine Steiner plane.
\end{example}

\begin{question} Is every finite balanced affine liner $X$ of rank $\|X\|\ge 4$ regular?
\end{question}

\begin{example}[{\cite[4.1.8]{Banakh}}] For every cardinal $\kappa$ there exists a non-regular $\kappa$-long Proclus liner $X$ of rank $\|X\|=4$. To construct such a liner, take any $(\kappa+1)$-long projective liner $Y$ of rank $\|Y\|=4$, choose a line $L\subseteq Y$ and a point $x\in L$. Then the subliner $X\defeq Y\setminus(L\setminus\{x\})$ is nonregular, $\kappa$-long, Proclus, and has rank $\|X\|=\|Y\|=4$.
\end{example}



\begin{theorem}[{\cite[4.5.1]{Banakh}}]\label{t:Playfair<=>regular} For a Playfair liner $X$ the following conditions are equivalent:
\begin{enumerate}
\item[\textup{(1)}] $X$ is regular;
\item[\textup{(2)}] $X$ is weakly regular;
\item[\textup{(3)}] $X$ is weakly modular;
\item[\textup{(4)}] $X$ is ranked;
\item[\textup{(5)}] $X$ is balanced;
\item[\textup{(6)}] $X$ is $4$-ranked or $4$-balanced;
\item[\textup{(7)}] $X$ contains no flat $F\subseteq X$ of rank $\|F\|=4$ and cardinality $|F|=81$.
\end{enumerate}
\end{theorem}

By Theorem~\ref{t:Playfair<=>}, a liner is Playfair if and only if it is affine, $3$-regular and $3$-long. In fact the first two conditions can be replaced by a single property, called the biaffinity.

\begin{definition} A liner $X$ is called \defterm{biaffine} if for every points $o,a,c\in X$, $b\in\Aline ac$ and $y\in\Aline ob\setminus\Aline ac$, there exist unique points $x\in\Aline oa$ and $z\in\Aline oc$ such that $y\in\Aline xz$ and $\Aline xz\cap\Aline ac=\varnothing$.
\end{definition}

\begin{picture}(200,80)(-200,-10)
\put(-30,0){\line(1,0){60}}
\put(-30,0){\line(1,2){30}}
\put(30,0){\line(-1,2){30}}
\put(-20,20){\line(1,0){40}}
\put(0,0){\line(0,1){60}}

\put(0,0){\circle*{3}}
\put(-3,-11){$b$}
\put(-30,0){\circle*{3}}
\put(-34,-9){$a$}
\put(30,0){\circle*{3}}
\put(29,-9){$c$}
\put(0,60){\circle*{3}}
\put(-2,63){$o$}
\put(-20,20){\color{red}\circle*{3}}
\put(-29,18){$x$}
\put(0,20){\circle*{3}}
\put(2,13){$y$}
\put(20,20){\color{red}\circle*{3}}
\put(23,18){$z$}
\end{picture}

\begin{theorem}[{\cite[4.5.3]{Banakh}}]\label{t:Playfair2<=>} A liner $X$ is Playfair if and only if $X$ is biaffine and $3$-long.
\end{theorem}

The following famous example of Marshall Hall \cite{Hall1960} shows that Theorem~\ref{t:4-long-affine} cannot be generalized to $3$-long affine liners.

\begin{example}[Hall, 1960; {\cite[4.1.11]{Banakh}}]\label{ex:HTS} There exists a Playfair liner, which is not $4$-ranked.
\end{example}

\begin{proof}[Proof] Let $\IF_3$ be the $3$-element field, $X\defeq\IF_3^4$ be the vector space over $\IF_3$ with origin $\boldsymbol o$ and the standard basis $\boldsymbol e_0,\boldsymbol e_1,\boldsymbol e_2,\boldsymbol e_3$. Let  $\mathcal L$ be the family of all 3-elements sets $\{\boldsymbol x,\boldsymbol y,\boldsymbol z\}\subseteq \IF_3^4$ such that $\boldsymbol x+\boldsymbol y+\boldsymbol z=(x_1-y_1)(x_2y_3-x_3y_2)\boldsymbol e_0$, where $(x_0,x_1,x_2,x_3)$ and $(y_0,y_1,y_2,y_3)$ are the coordinates of the vectors $\boldsymbol x,\boldsymbol y$ in the basis $\boldsymbol e_0,\boldsymbol e_1,\boldsymbol e_2,\boldsymbol e_3$. 
It can be shown that $(X,\mathcal L)$ is a non-regular Playfair liner. This liner contains pairwise disjoint lines ${\color{teal}L_1}=\boldsymbol e_1+\boldsymbol e_2+\IF_3\boldsymbol e_3$, ${\color{red}L_2}=2\boldsymbol e_1+2\boldsymbol e_2+\IF_3\boldsymbol e_3$,
${\color{blue}L_3}=2\boldsymbol e_1+\IF_3\boldsymbol e_3$ such that the set $\Aline{\color{teal}L_1}{\color{red}L_2}$ is a plane, equal to the union of three disjoint lines $\color{teal}L_1$, $\color{red}L_2$ and ${\color{purple}L_4}\defeq\IF_3(\boldsymbol e_0+\boldsymbol e_3)$. On the other hand, the set $\Aline {\color{teal}L_1}{\color{blue}L_3}$ is the union of the lines $\color{teal}L_1$, $\color{blue}L_2$ and the plane ${\color{cyan}P}\defeq\IF_3\boldsymbol e_0+2\boldsymbol e_2+\IF_3\boldsymbol e_3$. Consequently, $\overline{{\color{teal}L_1}\cup{\color{blue}L_3}}=X$ and $\|X\|=4$. On the other hand, $\{0\}\times\IF_3^3$ is a proper $3$-dimensional flat in $X$, witnessing that the liner $X$ is not $4$-ranked and hence is not regular. The plane $\overline{{\color{teal}L_1}\cup{\color{red}L_2}}$ has one-point intersection with the plane $\IF_3\boldsymbol e_2+\IF_3\boldsymbol e_4$, which is impossible in 3-dimensional regular spaces, see Theorem~\ref{t:w-modular<=>}.

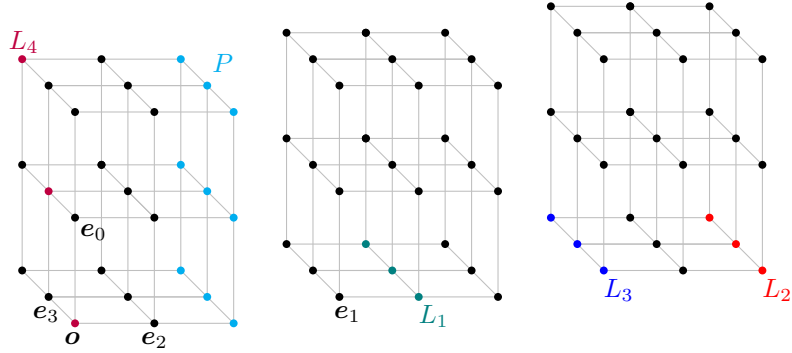
\begin{picture}(100,145)(-100,-65)


\put(-30,40){\color{lightgray}\line(1,0){60}}
\put(-30,-40){\color{lightgray}\line(1,0){60}}

\put(0,-40){\color{lightgray}\line(0,1){80}}
\put(20,-30){\color{lightgray}\line(0,1){80}}
\put(30,-40){\color{lightgray}\line(0,1){80}}
\put(40,-50){\color{lightgray}\line(0,1){80}}
\put(-20,-50){\color{lightgray}\line(0,1){80}}
\put(-30,-40){\color{lightgray}\line(0,1){80}}
\put(-40,-30){\color{lightgray}\line(0,1){80}}
\put(-10,-30){\color{lightgray}\line(0,1){80}}
\put(10,-50){\color{lightgray}\line(0,1){80}}

\put(-40,10){\color{lightgray}\line(1,0){60}}
\put(-30,0){\color{lightgray}\line(1,0){60}}
\put(-20,-10){\color{lightgray}\line(1,0){60}}
\put(-20,-10){\color{lightgray}\line(-1,1){20}}
\put(10,-10){\color{lightgray}\line(-1,1){20}}
\put(40,-10){\color{lightgray}\line(-1,1){20}}

\put(-40,50){\color{lightgray}\line(1,0){60}}
\put(-30,40){\color{lightgray}\line(1,0){60}}
\put(-20,30){\color{lightgray}\line(1,0){60}}
\put(-20,30){\color{lightgray}\line(-1,1){20}}
\put(10,30){\color{lightgray}\line(-1,1){20}}
\put(40,30){\color{lightgray}\line(-1,1){20}}
\put(10,-10){\color{lightgray}\line(-1,1){10}}

\put(-40,-30){\color{lightgray}\line(1,0){60}}
\put(-30,-40){\color{lightgray}\line(1,0){60}}
\put(-20,-50){\color{lightgray}\line(1,0){60}}
\put(-20,-50){\color{lightgray}\line(-1,1){20}}
\put(10,-50){\color{lightgray}\line(-1,1){20}}
\put(40,-50){\color{lightgray}\line(-1,1){20}}


\put(0,0){\circle*{3}}

\put(20,10){\color{cyan}\circle*{3}}
\put(30,0){\color{cyan}\circle*{3}}
\put(40,-10){\color{cyan}\circle*{3}}
\put(-10,10){\circle*{3}}
\put(-20,-10){\circle*{3}}
\put(-18,-17){$\boldsymbol e_0$}
\put(-30,0){\color{purple}\circle*{3}}
\put(-40,10){\circle*{3}}
\put(-10,10){\circle*{3}}
\put(10,-10){\circle*{3}}

\put(0,40){\circle*{3}}
\put(20,50){\color{cyan}\circle*{3}}
\put(32,45){\color{cyan}$P$}
\put(30,40){\color{cyan}\circle*{3}}
\put(40,30){\color{cyan}\circle*{3}}
\put(-10,50){\circle*{3}}

\put(-20,30){\circle*{3}}
\put(-30,40){\circle*{3}}
\put(-40,50){\color{purple}\circle*{3}}
\put(-45,54){\color{purple}$L_4$}
\put(10,30){\circle*{3}}

\put(0,-40){\circle*{3}}
\put(20,-30){\color{cyan}\circle*{3}}
\put(30,-40){\color{cyan}\circle*{3}}
\put(40,-50){\color{cyan}\circle*{3}}
\put(-20,-50){\color{purple}\circle*{3}}
\put(-24,-58){$\boldsymbol o$}
\put(-30,-40){\circle*{3}}
\put(-36,-48){$\boldsymbol e_3$}
\put(-40,-30){\circle*{3}}
\put(-10,-30){\circle*{3}}
\put(10,-50){\circle*{3}}
\put(5,-58){$\boldsymbol e_2$}


\put(70,50){\color{lightgray}\line(1,0){60}}
\put(70,-30){\color{lightgray}\line(1,0){60}}

\put(100,-30){\color{lightgray}\line(0,1){80}}
\put(120,-20){\color{lightgray}\line(0,1){80}}
\put(130,-30){\color{lightgray}\line(0,1){80}}
\put(140,-40){\color{lightgray}\line(0,1){80}}
\put(80,-40){\color{lightgray}\line(0,1){80}}
\put(70,-30){\color{lightgray}\line(0,1){80}}
\put(60,-20){\color{lightgray}\line(0,1){80}}
\put(90,-20){\color{lightgray}\line(0,1){80}}
\put(110,-40){\color{lightgray}\line(0,1){80}}

\put(60,20){\color{lightgray}\line(1,0){60}}
\put(70,10){\color{lightgray}\line(1,0){60}}
\put(80,0){\color{lightgray}\line(1,0){60}}
\put(80,0){\color{lightgray}\line(-1,1){20}}
\put(110,0){\color{lightgray}\line(-1,1){20}}
\put(140,0){\color{lightgray}\line(-1,1){20}}

\put(60,60){\color{lightgray}\line(1,0){60}}
\put(70,50){\color{lightgray}\line(1,0){60}}
\put(80,40){\color{lightgray}\line(1,0){60}}
\put(80,40){\color{lightgray}\line(-1,1){20}}
\put(110,40){\color{lightgray}\line(-1,1){20}}
\put(140,40){\color{lightgray}\line(-1,1){20}}

\put(60,-20){\color{lightgray}\line(1,0){60}}
\put(70,-30){\color{lightgray}\line(1,0){60}}
\put(80,-40){\color{lightgray}\line(1,0){60}}
\put(80,-40){\color{lightgray}\line(-1,1){20}}
\put(110,-40){\color{lightgray}\line(-1,1){20}}
\put(140,-40){\color{lightgray}\line(-1,1){20}}

\put(100,10){\circle*{3}}

\put(120,20){\circle*{3}}
\put(130,10){\circle*{3}}
\put(140,0){\circle*{3}}
\put(90,20){\circle*{3}}
\put(80,0){\circle*{3}}
\put(70,10){\circle*{3}}
\put(60,20){\circle*{3}}
\put(90,20){\circle*{3}}
\put(110,0){\circle*{3}}
\put(100,50){\circle*{3}}
\put(120,60){\circle*{3}}
\put(130,50){\circle*{3}}
\put(140,40){\circle*{3}}
\put(90,60){\circle*{3}}
\put(80,40){\circle*{3}}
\put(70,50){\circle*{3}}
\put(60,60){\circle*{3}}
\put(90,60){\circle*{3}}
\put(110,40){\circle*{3}}

\put(100,-30){\color{teal}\circle*{3}}

\put(120,-20){\circle*{3}}
\put(130,-30){\circle*{3}}
\put(140,-40){\circle*{3}}
\put(80,-40){\circle*{3}}
\put(78,-48){$\boldsymbol e_1$}
\put(70,-30){\circle*{3}}
\put(60,-20){\circle*{3}}
\put(90,-20){\color{teal}\circle*{3}}
\put(110,-40){\color{teal}\circle*{3}}
\put(110,-50){\color{teal}$L_1$}


\put(170,60){\color{lightgray}\line(1,0){60}}
\put(170,-20){\color{lightgray}\line(1,0){60}}

\put(200,-20){\color{lightgray}\line(0,1){80}}
\put(220,-10){\color{lightgray}\line(0,1){80}}
\put(230,-20){\color{lightgray}\line(0,1){80}}
\put(240,-30){\color{lightgray}\line(0,1){80}}
\put(180,-30){\color{lightgray}\line(0,1){80}}
\put(170,-20){\color{lightgray}\line(0,1){80}}
\put(160,-10){\color{lightgray}\line(0,1){80}}
\put(190,-10){\color{lightgray}\line(0,1){80}}
\put(210,-30){\color{lightgray}\line(0,1){80}}

\put(160,30){\color{lightgray}\line(1,0){60}}
\put(170,20){\color{lightgray}\line(1,0){60}}
\put(180,10){\color{lightgray}\line(1,0){60}}
\put(180,10){\color{lightgray}\line(-1,1){20}}
\put(210,10){\color{lightgray}\line(-1,1){20}}
\put(240,10){\color{lightgray}\line(-1,1){20}}

\put(160,70){\color{lightgray}\line(1,0){60}}
\put(170,60){\color{lightgray}\line(1,0){60}}
\put(180,50){\color{lightgray}\line(1,0){60}}
\put(180,50){\color{lightgray}\line(-1,1){20}}
\put(210,50){\color{lightgray}\line(-1,1){20}}
\put(240,50){\color{lightgray}\line(-1,1){20}}

\put(160,-10){\color{lightgray}\line(1,0){60}}
\put(170,-20){\color{lightgray}\line(1,0){60}}
\put(180,-30){\color{lightgray}\line(1,0){60}}
\put(180,-30){\color{lightgray}\line(-1,1){20}}
\put(210,-30){\color{lightgray}\line(-1,1){20}}
\put(240,-30){\color{lightgray}\line(-1,1){20}}

\put(200,20){\circle*{3}}

\put(220,30){\circle*{3}}
\put(230,20){\circle*{3}}
\put(240,10){\circle*{3}}
\put(190,30){\circle*{3}}
\put(180,10){\circle*{3}}
\put(170,20){\circle*{3}}
\put(160,30){\circle*{3}}
\put(190,30){\circle*{3}}
\put(210,10){\circle*{3}}

\put(200,60){\circle*{3}}
\put(220,70){\circle*{3}}
\put(230,60){\circle*{3}}
\put(240,50){\circle*{3}}
\put(190,70){\circle*{3}}
\put(180,50){\circle*{3}}
\put(170,60){\circle*{3}}
\put(160,70){\circle*{3}}
\put(190,70){\circle*{3}}
\put(210,50){\circle*{3}}

\put(200,-20){\circle*{3}}
\put(220,-10){\color{red}\circle*{3}}
\put(230,-20){\color{red}\circle*{3}}
\put(240,-30){\color{red}\circle*{3}}
\put(240,-40){\color{red}$L_2$}
\put(180,-30){\color{blue}\circle*{3}}
\put(180,-40){\color{blue}$L_3$}
\put(170,-20){\color{blue}\circle*{3}}
\put(160,-10){\color{blue}\circle*{3}}
\put(190,-10){\circle*{3}}
\put(210,-30){\circle*{3}}

\end{picture}
\end{proof}


\section{Hyperaffine, hyperbolic, and injective liners.}\label{s:hyperbolic}

In this section we present some (non-trivial) examples of finite hyperaffine, hyperbolic, and injective liners. 
  
\begin{example}[{\cite[3.8.1]{Banakh}}]\label{ex:hyperaffine-nonPlayfair} There exists a hyperaffine non-affine $2$-balanced $4$-parallel plane with $|X|_2=4$ and $|X|=|X|_3=25$. All its 50 lines are listed in the following table taken from \cite[1.34.1]{HCD}:
\begin{center}
{\tt 00000000111111122222223333344445555666778899aabbil\\
134567ce34578cd34568de468bh679f78ag79b9aabcddecejm\\
298dfbhkea6g9kf7c9afkg5cgfihdgifchi8ejjcjdfhgfghkn\\
iaolgmjnmbohnljonblhmjjdlknmeklnekmkinlimimonooloo
}
\end{center}
\end{example}


\begin{example}[{\cite[3.8.3]{Banakh}}]\label{ex:non-hyperaffine} There exists a $4$-parallel balanced ranked liner $X$ with $|X|_2=4$ and $|X|=|X|_3=25$, which is not hyperaffine. It is the group $\IZ_5\times\IZ_5$ endowed with the family of lines $\mathcal L\defeq\{B+z:B\in\big\{\{00,10,01,22\},\{13,31,33,44\},\;z\in \IZ_5\times \IZ_5\}$. The non-hyperaffinity of this liner is witnessed by the points $o\defeq 00$, $x\defeq 41$, $y\defeq 01$,  $p\defeq 13$.
\end{example}

\begin{example}[{\cite[3.8.4]{Banakh}}]\label{ex:Hetman-hyperbolic91} There exists a hyperbolic $8$-parallel plane $X$ with $|X|_2=7$ and $|X|=|X|_3=91$. It is the ring $\IZ_{91}$ endowed with the family of lines $\mathcal L\defeq\{x+L:\break L\in\big\{\{0, 8, 29, 51, 54, 61, 63\},\{0, 11, 16, 17, 31, 35, 58\},\{0, 13, 26, 39, 52, 65, 78\}\big\},\;x\in\IZ_{91}\}.$
\end{example}

\begin{example}[{\cite[3.8.7]{Banakh}}]\label{ex:Hetman-hyperbolic} For every $n\in\{3,4,5,7,8,9\}$, there exists an injective $2$-balanced $(n^2-n-1)$-parallel plane $X$ with $|X|_2=n+1$ and $|X|=|X|_3=n^3+1$. 
\end{example}

The injective (and hence hyperbolic) liners appearing in Example~\ref{ex:Hetman-hyperbolic} are classical unitals.

\begin{definition} A \defterm{unital} is a $2$-balanced liner $X$ with $|X|=n^3+1$ and $|X|_2=n+1$ for some finite cardinal $n\ge 2$. 
\end{definition}


\begin{definition} A \defterm{classical unital} is the subliner $$U_q\defeq \{\IF_{q^2}(x,y,z):(x,y,z)\in\IF_{q^2}^3\setminus\{0\}^3\;\;(x^{q+1}+y^{q+1}+z^{q+1}=0)\}$$ 
of the projective plane $\mathbb P\IF_{q^2}^3$ over a field $\IF_{q^2}$ whose order $|\IF_{q^2}|=q^2$ is a square of some number $q$ (which is necessarily a prime power).
\end{definition}

\begin{remark}[{\cite[7.42]{UPP}}] Every classical unital $U_q$ is an injective $2$-balanced $(q^2-q-1)$-parallel liner with $|U_q|_2=q+1$ and $|U_q|=|U_q|_3=q^3+1$. If $q=2$, then the liner $U_q$ is affine. If $q\ge 3$, then the liner $U_q$ is hyperbolic.
\end{remark} 

In fact, injective liners are not exotic at all, and as subliners are present in every $\w$-long plane. We recall that a plane is a liner of rank $3$.

\begin{theorem}[{\cite[3.8.11]{Banakh}}]\label{t:injective-subliner} Let $X$ be an $\w$-long plane. For every cardinal $\lambda\in [3,\w]$, there exists an injective subliner $Y\subseteq X$ such that 
\begin{enumerate}
\item[\textup{(1)}] $|Y|=\|Y\|=\w$;
\item[\textup{(2)}] $Y$ is $n$-balanced for every positive cardinal $n$;
\item[\textup{(3)}] $|Y|_2=\lambda$ and $|Y|_n=\w$ for every cardinal $n\ge 3$.
\end{enumerate}
\end{theorem}

\begin{question}\label{q:hyperbolic} Is every $2$-hyperparallel {\em regular} liner $X$ hyperaffine? hyperbolic?
\end{question}

\begin{remark} Examples~\ref{ex:Steiner13} and \ref{ex:non-hyperaffine} show that the regularity is essential in Question~\ref{q:hyperbolic}.
\end{remark}

\section{Steiner liners}

By Theorem~\ref{t:4-long-affine}, every $4$-long affine liner is regular. If an affine liner $X$ is not $4$-long, then all lines in $X$ have the same length $|X|_2\in\{2,3\}$. 
Finite $2$-balanced liners $X$ with $|X|_2=3$ are well-known in Combinatorics as Steiner Triple Systems, see \cite[\S II.2]{HCD}. This motivates the following definition.

\begin{definition} A liner $X$ is called \defterm{Steiner} if every line in $X$ contains exactly three points.
\end{definition}

It is easy to see that all Steiner liners are proaffine. Moreover, a Steiner liner is affine if and only if it is hyperaffine if and only if it is injective.

Theorem~\ref{t:injective-subliner} implies the following corollary supplying us with many examples of infinite affine Steiner liners.

\begin{corollary}[{\cite[4.2.4]{Banakh}}] Let $r\in\{3,\w\}$ Every $\w$-long plane $P$ contains an infinite affine Steiner subliner $X$ of rank $\|X\|=r$. 
\end{corollary}

On the other hand, Theorem~\ref{t:2-balanced-dependence} implies that the following corollary imposing restrictions on possible cardinalities of finite Steiner liners.

\begin{corollary}[Kirkman, 1846; {\cite[4.2.5]{Banakh}}]\label{c:Kirkman} Every finite Steiner liner $X$ has cardinality\newline $|X|\in (6\IN+1)\cup(6\IN-3).$
\end{corollary} 

\begin{remark}The necessary condition given by Corollary~\ref{c:Kirkman} is also sufficient: for every number $n\in (6\IN+1)\cup(6\IN-3)$ there exists a Steiner liner of cardinality $n$, see \cite{Bose1939}, \cite{Skolem1958}, \cite{RCW1971}. Moreover, if $n\in (6\IN+13)\cup(6\IN-3)$, then there exists an {\em affine} Steiner liner of cardinality $n$, see \cite{GGW}.
\end{remark}

\begin{theorem}[{\cite[4.2.7]{Banakh}}]\label{t:Steiner-3reg} Every Steiner $3$-regular liner $X$ is $p$-parallel for some $p\in\{0,1\}$. Consequently, $X$ is projective or affine.
\end{theorem}

 Every Steiner liner $X$ carries the commutative binary operation $\circ:X\times X\to X$ assigning to every pair $(x,z)\in X\times X$ the unique point $y\in X$ such that $\{x,y,z\}=\Aline xz$. This operation will be called the \defterm{midpoint operation} on $X$. The midpoint operation satisfies the identities 
$$ x\circ x=x,\quad x\circ y=y\circ x,\quad (x\circ y)\circ y=x,$$
turning $X$ into an involutory idempotent commutative magma. 
\smallskip

Now we recall the definitions of some algebraic structures, starting with the most basic notion of a magma. A \defterm{magma} is a set $X$ endowed with a binary operation $\cdot:X\times X\to X$, $\cdot:(x,y)\mapsto xy$.

A magma $X$ is 
\begin{itemize}
\item \defterm{commutative} if $xy=yx$ for every $x,y\in X$;
\item \defterm{associative} if $(xy)z=x(yz)$ for every $x,y,z\in X$;
\item \defterm{idempotent} if $xx=x$ for every $x\in X$;
\item \defterm{involutory} if $x(xy)=y$ for every $x,y\in X$;
\item \defterm{self-distributive} if $x(yz)=(xy)(xz)$;
\item \defterm{unital} if there exists an element $e\in X$ such that $xe=x=ex$ for every $x\in X$;
\item a \defterm{quasigroup} if for every $a,b\in X$ there exist unique $x,y\in X$ such that $ax=b=ya$;
\item a \defterm{loop} if $X$ is a unital quasigroup;
\item a \defterm{group} of $X$ is an associative loop;
\item \defterm{Moufang} if $(xy)(zx)=(x(yz))x=x((yz)x)$ for every $x,y\in X$.
\end{itemize}

A subset $S$ of a magma $X$ is called a \defterm{submagma} of the magma $X$ if $\{xy:x,y\in S\}\subseteq S$. A subset $S$ of a loop $X$ is called a \index{subloop}\defterm{subloop} of $X$ if for every $a,b\in S$, $ab\in S$ and there exist points $x,y\in S$ such that $ax=b=ya$.

\begin{theorem}[{\cite[4.2.8]{Banakh}}]\label{t:Steiner<=>midpoint} Any Steiner liner endowed with the midpoint operation is an involutory idempotent com\-mu\-ta\-tive magma. Conversely, every involutory idempotent commutative magma $(M,\circ)$ endowed with the family of lines $\mathcal L\defeq\{\{x,y,x\circ y\}:x,y\in M\;\wedge\;x\ne y\}$ is a Steiner liner. Moreover, a subset $F\subseteq X$ is flat in the liner $X$ iff $F$ is a submagma of the magma $(X,\circ)$.
\end{theorem}

A less evident is the relation of Steiner liners to commutative loops. 
By definition, every loop $X$ contains an element $e\in X$ such that $ex=x=xe$ for every $x\in X$. This element is unique and is called the \defterm{identity element} of the loop. A loop $X$ is defined to be of \defterm{exponent $3$} if  $(xx)x=x(xx)=e$  for every element $x\in X$.

For every element $e$ of a Steiner liner $X$, the binary operation
$$\cdot:X\times X\to X,\quad \cdot:(x,y)\mapsto xy\defeq e\circ(x\circ y),$$
turns $X$ into a commutative loop of exponent $3$ with identity element $e$. 

We shall denote this loop by $X_e$. 

\begin{proposition}[{\cite[4.2.10]{Banakh}}]\label{p:flat<=>subloop} For any set $F$ in a Steiner liner $X$ and any element $e\in F$, the following conditions are equivalent:
\begin{enumerate}\setlength{\itemsep}{-0pt}\setlength{\parskip}{1pt}
\item[\textup{(1)}] $F$ is a flat in the liner $X$;
\item[\textup{(2)}] $F$ is a subloop of the loop $X_e$.
\item[\textup{(3)}] $F$ is a submagma of the magma $X_e$;
\item[\textup{(4)}] $F$ is a submagma of the magma $(X,\circ)$.
\end{enumerate}
\end{proposition}

\begin{theorem}[{\cite[4.2.12]{Banakh}}]\label{t:Steiner<=>projective} For a Steiner liner $X$ the following conditions are equivalent:
\begin{enumerate}\setlength{\itemsep}{-0pt}\setlength{\parskip}{1pt}
\item[\textup{(1)}] $X$ is projective;
\item[\textup{(2)}]  $X$ is $0$-parallel;
\item[\textup{(3)}]  $X$ is $3$-balanced with $|X|_3=7$;
\item[\textup{(4)}]  for any distinct points $x,y,z\in X$ we have $(x\circ y)\circ (x\circ z)=y\circ z$.
\end{enumerate}
\end{theorem}

\section{Hall liners}\label{s:Hall}

\begin{definition} A liner $X$ is called \defterm{Hall}  if $X$ is $2$-balanced and $3$-balanced with
\smallskip

\centerline{$|X|_3=(|X|_2)^2=9.$}
\end{definition}

The liner from Example~\ref{ex:HTS} is an example of a Hall liner. Hall liners are well-known in Theory of Combinatorial Designs as \defterm{Hall Triple Systems}, see \cite[\S VI.28]{HCD}. Corollary~\ref{c:Playfair<=>23balanced} implies the following characterization of Hall liners.

\begin{theorem}[{\cite[4.4.2]{Banakh}}]\label{t:Hall<=>Playfair+Steiner} A liner $X$ is Hall if and only if $X$ is Steiner and Playfair.
\end{theorem}

\begin{theorem}[{\cite[4.4.3]{Banakh}}]\label{t:Hall<=>Moufang} For a Steiner liner $X$, the following conditions are equivalent:
\begin{enumerate}\setlength{\itemsep}{-0pt}\setlength{\parskip}{1pt}
\item[\textup{(1)}] $X$ is Hall;
\item[\textup{(2)}]  the mid-point operation on $X$ is self-distributive;
\item[\textup{(3)}]  for every point $e\in X$, the commutative loop $X_e$ is Moufang.
\end{enumerate}
\end{theorem}

\begin{theorem}[{\cite[4.4.5]{Banakh}}]\label{t:Hall<=>regular} For a Hall liner $X$, the following conditions are equivalent:
\begin{enumerate}\setlength{\itemsep}{-0pt}\setlength{\parskip}{1pt}
\item[\textup{(1)}] $X$ is regular;
\item[\textup{(2)}] $X$ is weakly regular;
\item[\textup{(3)}] $X$ is weakly modular;
\item[\textup{(4)}] $X$ is $4$-regular;
\item[\textup{(5)}] $X$ is $4$-ranked;
\item[\textup{(6)}] $X$ is $4$-balanced;
\item[\textup{(7)}] $X$ is $4$-balanced with $|X|_4=27$;
\item[\textup{(8)}] $X$ contains no flat $F\subseteq X$ of rank $\|F\|=4$ and cardinality $|F|=81$;
\item[\textup{(9)}] for every point $e\in X$, the loop $X_e$ is associative.
\end{enumerate}
\end{theorem}

A group $X$ is defined to be \defterm{nilpotent of class $n$} for $n\in\IN$ if the quotient group $X/Z$ of $X$ by its centre $Z\defeq\{z\in X:\forall x\in X\;(xz=zx)\}$ is nilpotent of class  $n-1$. The trivial group is nilpotent of class $0$. So, non-trivial commutative groups are nilpotent of class $1$. A group $X$ is nilpotent of class at most $2$ if and only if the quotient group $X/Z$ of $X$ by its centre $Z$ is commutative.   
It is known \cite{Burnside1901}, \cite{Levi1942} that every group $X$ of exponent $3$ is nilpotent of class at most $3$. In particular, for every $m\in \IN$, the free Burnside group $B(m,3)=\langle x_1,\dots,x_m:x^3=1\rangle$ of exponent $3$ with $m$ generators is nilpotent of class at most $3$. The Burnside group $B(1,3)$ is nilpotent of class $1$ and has cardinality $3$, the Burnside group $B(2,3)$ is nilpotent of class 2 with $|B(2,3)|=3^3=27$, and for every $m\ge 3$, the Burnside group $B(m,3)$ is nilpotent of class $3$ with $|B(m,3)|=3^{C^1_m+C^2_m+C^3_m}$.

\begin{example}[{\cite[4.4.7]{Banakh}}] Every group $X$ of exponent $3$ endowed with the binary operation
\smallskip

\centerline{$\circ:X\times X\to X,\quad\circ:(x,y)\mapsto x\circ y\defeq xy^{-1}x,$} 
\smallskip

\noindent is a self-distributive idempotent involutory commutative magma. Consequently, $X$ endowed with the family of lines $\mathcal L\defeq\big\{\{\{x,y,xy^{-1}x\}:x,y\in X\;\wedge\;x\ne y\big\}$ is a Hall liner. This Hall liner is regular if and only if the group $X$ is nilpotent of class at most $2$. 
\end{example}

\begin{remark} By a deep result of Bruck and Slaby \cite[Theorem 10.1]{Bruck1958}, every finitely generated commutative Moufang loop $X$ of exponent $3$ is centrally nilpotent, which implies that $|X|=3^n$ for some $n\in\IN$. This algebraic result implies that every Hall liner $X$ of finite rank has cardinality $|X|=3^n$ for some $n\in\IN$. This fact was also proved by Hall \cite{Hall1980} by group-theoretic methods, without involving the machinery of Moufang loops. Theorem~\ref{t:Hall<=>regular} implies that every non-regular Hall liners of rank 4 are isomorphic and have cardinality $81$. The following table taken from \cite[\S VI.28.10]{HCD} shows the number $\hbar$ of non-isomorphic non-regular Halls liners of cardinality $3^n$.
\begin{center}
\begin{tabular}{c||c|c|c|c|c|c|c|c|c}
$n$&1&2&3&4&5&6&7&8&9\\
\hline
$\hbar$&0&0&0&1&1&3&12&$\ge 45$&?
\end{tabular}
\end{center}
The liner in Example~\ref{ex:HTS} is the unique non-regular Hall liner of cardinality $81$.
\end{remark}

\section{Subparallelity}\label{ss:subparallelity-in-matroids}

\begin{definition}
  Given two flats $A$ and $B$ in a liner, we write $A\subparallel B$ and say that $A$ is \defterm{subparallel to} $B$ if $\forall a\in A\;\;(A\subseteq \overline{\{a\}\cup B})$.
\end{definition}

\begin{lemma}[{\cite[6.1.2]{Banakh}}]\label{l:subparallel+intersect=>subset} For two intersecting flats $A,B$ be in a liner, we have $A\subparallel B$ if and only if $A\subseteq B$.
\end{lemma}

It can be shown that a flat $A$ in a liner $X$ is subparallel to a flat $B\subseteq X$ if and only if for every $a,b\in A$ the flat $\Aline ab$ is subparallel to $B$.

\begin{theorem}[{\cite[6.1.4]{Banakh}}]\label{t:subparallel-char} Let $X$ be a $\kappa$-ranked liner for some cardinal $\kappa$. For two flats $A,B\subseteq X$ with $\|B\|<\kappa$, the following conditions are equivalent:
\begin{enumerate}\setlength{\itemsep}{-0pt}\setlength{\parskip}{1pt}
\item[\textup{(1)}] $A\subparallel B$;
\item[\textup{(2)}] either $A\subseteq B$ or  $A\subseteq \overline{\{a\}\cup B}\setminus B$ for some $a\in A$;
\item[\textup{(3)}] either $A\subseteq B$ or  $A\cap B=\varnothing$ and $\|A\|_B=1$.
\end{enumerate}
If $\|B\|<\w$, then the conditions \textup{(1)--(3)} are equivalent to the condition
\begin{itemize}
\item[\textup{(4)}] either $A\subseteq B$ or $A\cap B=\varnothing$ and $\|A\cup B\|=\|B\|+1$.
\end{itemize}
\end{theorem}

\begin{corollary}[{\cite[6.1.5]{Banakh}}]\label{c:subparallel} Let $\kappa$ be a finite cardinal and $A,B$ are two flats in a $\kappa$-ranked liner $X$. If $\|A\|=\|B\|<\kappa$, then $A\subparallel B\;\Leftrightarrow\; B\subparallel A$.
\end{corollary}

\begin{corollary}[{\cite[6.1.7]{Banakh}}]\label{c:para+intersect=>coincide} Let $A,B$ be two flats in a $\|B\|$-ranked liner $X$ such that $A\subparallel B$ and $\|A\|=\|B\|<\w$. If $A\cap B\ne\varnothing$, then $A=B$.
\end{corollary}

\begin{proposition}[{\cite[6.1.8]{Banakh}}]\label{p:subparallel=>dim} Let $\kappa$ be a cardinal and $A,B$ be two flats in a $\kappa$-ranked liner $X$. If $\|B\|<\kappa$, $A\subparallel B$ and $B\ne\varnothing$, then $\|A\|\le \|B\|$.
\end{proposition}

\begin{remark} The non-ranked liner in Example~\ref{ex:Tao} contains two lines $L,\Lambda$ such that $L\subparallel \Lambda$ but not $\Lambda\subparallel L$.
\end{remark}

\section{Parallelity}\label{ss:parallelity-in-matroids}

\begin{definition} Given two flats $A$ and $B$ in a liner $X$, we write $A{\parallel}B$ and say that $A$ and $B$ are \defterm{parallel} if $A$ is subparallel to $B$  and $B$ is subparallel to $A$. Therefore,
$A{\parallel}B\;\Leftrightarrow\;(A\subparallel B\;\wedge\; B\subparallel A).$ For two flats $A,B$, the negation of $A\parallel B$ is denoted by $A\nparallel B$.
\end{definition}

Theorem~\ref{t:subparallel-char} implies the following characterization.

\begin{theorem}[{\cite[6.2.2]{Banakh}}]\label{t:parallel-char} Let $\kappa$ be a cardinal, $X$ be a $\kappa$-ranked liner, and $A,B$ be two nonempty flats in $X$. If $\max\{\|A\|,\|B\|\}<\kappa$, then the following conditions are equivalent:
\begin{enumerate}\setlength{\itemsep}{-0pt}\setlength{\parskip}{1pt}
\item[\textup{(1)}] $A\parallel B$;
\item[\textup{(2)}] either $A=B$ or  $A\subseteq \overline{\{a\}\cup B}\setminus B$ and $B\subseteq\overline{\{b\}\cup A}\setminus A$ for some points $a\in A$, $b\in B$;
\item[\textup{(3)}] either $A=B$ or else $A\cap B=\varnothing$ and $\|A\|_B=1=\|B\|_A$.
\end{enumerate}
If $\min\{\|A\|,\|B\|\}<\w$, then the conditions \textup{(1)--(3)} are equivalent to the condition
\begin{enumerate}
\item[\textup{(4)}] either $A=B$ or else $A\cap B=\varnothing$ and $1+\|A\|=\|A\cup B\|=1+\|B\|$.
\end{enumerate}
\end{theorem}

\begin{remark} Any two singletons in a liner are parallel.
\end{remark}

Proposition~\ref{p:subparallel=>dim} implies

\begin{corollary}[{\cite[6.2.4]{Banakh}}]\label{c:parallel} Let $\kappa$ be a cardinal, $X$ be a $\kappa$-ranked liner and $A,B$ be two flats in $X$. If $0<\min\{\|A\|,\|B\|\}<\kappa$ and $A\parallel B$, then $\|A\|=\|B\|$.
\end{corollary}

Lemma~\ref{l:subparallel+intersect=>subset} implies the following simple (but helpful) proposition.

\begin{proposition}[{\cite[6.2.5]{Banakh}}]\label{p:para+intersect=>coincide} Two intersecting flats in a liner are parallel if and only if they are equal.
\end{proposition}

Theorem~\ref{t:parallel-char} implies the following characterization of the parallelity of lines in $3$-ranked liner, which often is taken as the definition of the parallelity of lines.

\begin{corollary}[{\cite[6.2.6]{Banakh}}]\label{c:parallel-lines<=>} Two lines $L,\Lambda$ is a $3$-ranked liner are parallel if and only if $\|L\cup\Lambda\|\le 3$ and either $L=\Lambda$ or $L\cap\Lambda=\varnothing$.
\end{corollary}

\begin{corollary}[{\cite[6.2.7]{Banakh}}]\label{c:parallel-in-Proclus} Let $A,B,C$ be three lines in a Proclus liner $X$. If $A\parallel B$, $B\parallel C$ and $\|A\cup C\|\le 3$, then $A\parallel C$.
\end{corollary}

\section{Parallelity in (weakly) modular liners}

\begin{proposition}[{\cite[6.3.1]{Banakh}}] For two flats $A,B$ in a modular liner $X$, the following conditions are equivalent:
\begin{enumerate}\setlength{\itemsep}{-0pt}\setlength{\parskip}{1pt}
\item[\textup{(1)}] $A\subparallel B$;
\item[\textup{(2)}]  $A\subseteq B$ or $\|A\|=1$.
\item[\textup{(3)}]  $A\subseteq B$ or $\|A\|\le 1$.
\end{enumerate}
\end{proposition}

\begin{corollary}[{\cite[6.3.2]{Banakh}}]\label{c:parallel-in-modular} For two flats $A,B$ in a modular liner $X$, the following conditions are equivalent:
\begin{enumerate}\setlength{\itemsep}{-0pt}\setlength{\parskip}{1pt}
\item[\textup{(1)}] $A\parallel B$;
\item[\textup{(2)}]  $A=B$ or $\max\{\|A\|,\|B\|\}=1$;
\item[\textup{(3)}]  $A=B$ or $\max\{\|A\|,\|B\|\}\le 1$.
\end{enumerate}
\end{corollary}

\begin{proposition}[{\cite[6.4.1]{Banakh}}] Let $A,B,C$ be flats in a weakly modular liner $X$. 
 If $A\subparallel B$ and $B\cap C\ne\varnothing$, then $A\cap C\subparallel B\cap C$.
 \end{proposition}
 
\begin{corollary}[{\cite[6.4.2]{Banakh}}]\label{c:paraintersect} Let $A,B,C$ be flats in a weakly modular liner $X$. 
 If $A\parallel B$ and $A\cap C\ne\varnothing\ne B\cap C$, then $A\cap C\parallel B\cap C$.
\end{corollary}
 
\begin{proposition}[{\cite[6.4.3]{Banakh}}]\label{p:subparallel-char4} For two disjoint flats $A,B$ in a weakly modular liner $X$, the following conditions are equivalent:
\begin{enumerate}\setlength{\itemsep}{-0pt}\setlength{\parskip}{1pt}
\item[\textup{(1)}] $A\subparallel B$;
\item[\textup{(2)}] for every $b\in B$, the flats $A$ and $B\cap\overline{\{b\}\cup A}$ are parallel;
\item[\textup{(3)}] there exists a flat $C\subseteq B$ such that $A\parallel C$;
\item[\textup{(4)}] $\forall a,x\in A\;\forall b\in B\;\exists y\in B\;\;(\Aline ax\parallel \Aline by)$.
\end{enumerate}
\end{proposition}

\section{Parallelity in proaffine regular liners}

\begin{theorem}[{\cite[6.5.1]{Banakh}}]\label{t:proaffine-char} For a regular liner $X$, the following conditions are equivalent:
\begin{enumerate}\setlength{\itemsep}{-0pt}\setlength{\parskip}{1pt}
\item[\textup{(1)}] the liner $X$ is proaffine;
\item[\textup{(2)}]  for every flat $A\subseteq X$ and points $a\in A$, $b\in X\setminus A$, $p\in \overline{A\cup\{b\}}\setminus A$,\newline there exists a point $u\in\Aline ab$ such that $\forall v\in\Aline ab\setminus\{u\}\;\;(\Aline vp\cap A\ne\varnothing)$;
\item[\textup{(3)}]  for every flat $A\subseteq X$ and point $b\in X\setminus A$ there exists at most one flat $B\subseteq X$\newline such that $b\in B$ and $B\parallel A$.
\end{enumerate}
\end{theorem}

In proaffine regular liners, the parrallelism is an equivalence relation (which is a generalization of  Proposition I.30 on parallel lines in Euclid's ``Elements'' \cite{Euclid}).
 
\begin{theorem}[{\cite[6.5.4+6.5.6]{Banakh}}]\label{c:subparallel-transitive} For any flats $A,B,C$ in a proaffine regular liner $X$ we have the implications $(A{\subparallel}B\wedge B{\subparallel}C)\Ra (A\subparallel C)$ and $(A{\parallel}B\wedge B{\parallel}C)\Ra (A\parallel C).$ 
\end{theorem}

\begin{example} The Hall liner in Example~\ref{ex:HTS} contains lines $L_1,L_2,L_3$ such that $L_1\parallel L_2\parallel L_3$ but $L_1\nparallel L_3$.
\end{example}

The following theorem generalizes the Proclus Parallelity Postulate~\ref{p:Proclus-Postulate}.

\begin{theorem}[{\cite[6.5.8]{Banakh}}]\label{t:Proclus2} Let $A,B$ be parallel flats in a proaffine regular liner $X$ and $L$ be a line in the flat $\overline{A\cup B}$. If $|L\cap A|=1$, then $|L\cap B|=1$.
\end{theorem}

\begin{theorem}[{\cite[6.5.9]{Banakh}}]\label{t:subparallel-via-base} A nonempty flat $A$ in a proaffine regular liner $X$ is subparallel to a flat $B$ in $X$ if and only if there exists a point $a\in A $ and a set $\Lambda\subseteq A$ such that $A=\overline\Lambda$ and $\forall x\in \Lambda\;(\Aline ax\subparallel B)$.
\end{theorem}

\section{Parallelity in affine and Playfair liners}

\begin{theorem}[{\cite[6.6.1]{Banakh}}]\label{t:affine-char2} For a regular liner $X$, the following conditions are equivalent:
\begin{enumerate}\setlength{\itemsep}{-0pt}\setlength{\parskip}{1pt}
\item[\textup{(1)}] the liner $X$ is affine;
\item[\textup{(2)}] for every flat $A\subseteq X$ and points $a\in A$, $b\in X\setminus A$, $p\in \overline{A\cup\{b\}}\setminus A$,\newline there exists a unique point $u\in\Aline ab$ such that $\Aline up\cap A=\varnothing$;
\item[\textup{(3)}] for every flat $A\subseteq X$ and points $a\in A$, $b\in X\setminus A$ with $\Aline ab\ne\{a,b\}$,\newline there exists a unique flat $B\subseteq X$ such that $b\in B$ and $B\parallel A$.
\end{enumerate}
\end{theorem}

\begin{theorem}[{\cite[6.6.2]{Banakh}}]\label{t:Playfair} For every line $L$ in a Playfair liner $X$ and every point $x\in X\setminus L$, there exists a unique line $\Lambda$ in $X$ such that  $x\in\Lambda$ and $\Lambda\parallel L$.
\end{theorem}

For regular Playfair liners, we have a higher dimensional counterpart of Theorem~\ref{t:Playfair}, which follows from Theorem~\ref{t:affine-char2}.

\begin{corollary}[{\cite[6.6.3]{Banakh}}]\label{c:Playfair} For any flat $A$ in a regular Playfair liner $X$ and every point $x\in X$, there exists a unique flat $B\subseteq X$ such that $x\in B$ and $B\parallel A$.
\end{corollary}

Let us recall that a liner $X$ is \defterm{Bolyai} if for every plane $P\subseteq X$, line $L\subseteq P$ and point $x\in P\setminus L$, there exists a line $\Lambda$ in $X$ such that $x\in \Lambda\subseteq P\setminus L$.

\begin{theorem}[{\cite[6.6.4]{Banakh}}]\label{t:Play-reg<=>Bolyai} For every liner $X$ the following conditions are equivalent:
\begin{enumerate}\setlength{\itemsep}{-0pt}\setlength{\parskip}{1pt}
\item[\textup{(1)}] $X$ is Playfair and regular;
\item[\textup{(2)}] $X$ is Bolyai, $3$-ranked, and for every lines $A,B,C\subseteq X$,\newline if $A\parallel B$ and $B\parallel C$, then $A\parallel C$.
\end{enumerate}
\end{theorem}

\begin{proposition}[{\cite[6.6.5]{Banakh}}]\label{p:ABC-regular} For any parallel lines $A,B,C$ in a Playfair liner $X$, the flat $\overline{A\cup B\cup C}$ is a regular subliner of $X$.
\end{proposition}

\section{Triangles and parallelograms in liners}

\begin{definition} A triple of points $abc\in X^3$ in a liner $X$ is called a \defterm{triangle} in $X$ if $\|\{a,b,c\}\|=3$. The points $a,b,c$ are called the \defterm{vertices} of the triangle $abc$, and the lines $\Aline ab,\Aline bc,\Aline ca$ are called the \defterm{sides} of the triangle $abc$.
\end{definition}

\begin{definition}
A quadruple of points $abcd\in X^4$  in a liner $X$ is called a \defterm{parallelogram} in $X$ if $$\Aline ab\parallel\Aline cd\ne \Aline ab\quad\mbox{and}\quad\Aline bc\parallel\Aline ad\ne \Aline bc.$$ The points $a,b,c,d$ are called the \defterm{vertices}, the lines $\Aline ab,\Aline bc,\Aline cd,\Aline da$ are called the \defterm{sides}, and the lines $\Aline ac,\Aline bd$ are called the \defterm{diagonals} of the parallelogram $abcd$.

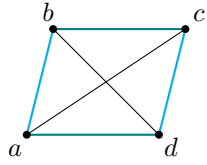
\begin{picture}(300,68)(-150,10)
{\linethickness{0.8pt}
\put(15,20){\color{teal}\line(1,0){50}}
\put(25,60){\color{teal}\line(1,0){50}}
\put(15,20){\color{cyan}\line(1,4){10}}
\put(65,20){\color{cyan}\line(1,4){10}}
}
\put(15,20){\line(3,2){60}}
\put(65,20){\line(-1,1){40}}
\put(15,20){\circle*{3}}
\put(8,12){$a$}
\put(65,20){\circle*{3}}
\put(67,12){$d$}
\put(25,60){\circle*{3}}
\put(21,63){$b$}
\put(75,60){\circle*{3}}
\put(78,63){$c$}
\end{picture}
\end{definition}

Corollary~\ref{c:parallel-lines<=>} implies the following useful characterization of parallelograms in $3$-ranked liners.

\begin{proposition}[{\cite[6.7.4]{Banakh}}]\label{p:tobe-parallelogram} A quadruple of points $abcd$ in a $3$-ranked liner $X$ is\newline a parallelogram if and only if  $\Aline ab\cap\Aline cd=\varnothing=\Aline bc\cap\Aline ad$ and $\|\{a,b,c,d\}\|\le 3$.
\end{proposition}

\begin{proposition}[{\cite[6.7.7]{Banakh}}]\label{p:k-parallel=>parallelogram} If a liner $X$ is  $\kappa$-parallel for some finite nonzero cardinal $\kappa$, then for every triangle $abc$ there exists a point $d\in X$ such that $abcd$ is a parallelogram.
\end{proposition}

\begin{theorem}[{\cite[6.7.8]{Banakh}}]\label{t:parallelogram3+1}  A liner $X$ is Playfair if and only if $X$ is affine, $3$-ranked, and for any triangle $abc$ in $X$ there exists a unique point $d\in X$ such that $abcd$ is parallelogram.
\end{theorem}

\section{Boolean parallelograms and Boolean liners}

\begin{definition}A parallelogram $abcd$ in a liner $X$ is called \defterm{Boolean} if its diagonals $\Aline ac$ and $\Aline bd$ are parallel. 
\end{definition}

For $3$-ranked Steiner affine liners, Proposition~\ref{p:k-parallel=>parallelogram} and Theorem~\ref{t:parallelogram3+1} are completed by the following proposition.

\begin{proposition}[{\cite[6.8.2]{Banakh}}]\label{p:Steiner-paragram} For any triangle $abc$ in a $3$-ranked Steiner affine liner $X$, there exists a unique point $d\in X$ such that $abcd$ is a non-Boolean parallelogram.
\end{proposition}

\begin{corollary}[{\cite[6.8.3]{Banakh}}] For every triangle $abc$ in a $3$-long $3$-ranked affine liner $X$, there exists a point $d\in X$ such that $abcd$ is a parallelogram.
\end{corollary} 
 
\begin{corollary}[{\cite[6.8.4]{Banakh}}]\label{c:Steiner-Playfair-Boolean} A Steiner liner is Playfair if and only if it is affine, $3$-ranked, and contains no Boolean parallelograms.
\end{corollary}

Let us recall that a liner is \defterm{line-finite} if all its lines are finite sets. Observe that a $2$-balanced liner $X$ is line-finite if and only if the cardinal $|X|_2$ is finite.

\begin{proposition}[{\cite[6.8.5]{Banakh}}]\label{p:large=>Boolean-parallelogram} If a line-finite $2$-balanced $3$-balanced $3$-ranked liner $X$ has $|X|_3\ge |X|_2^2$ \textup{(}and $|X|_3>3{\cdot}|X|_2^2-9{\cdot}|X|_2+9$\textup{)}, then for every triangle $abc$ there exists a point $d\in X$ such that $abcd$ is a \textup{(}Boolean\textup{)} parallelogram.
\end{proposition}

\begin{corollary}[{\cite[6.8.6]{Banakh}}] For every triangle $abc$ in a Steiner ranked plane $X$ of cardinality $|X|>9$, there exists a point $d\in X$ such that $abcd$ is a Boolean parallelogram.
\end{corollary}

\begin{definition} A liner $X$ is defined to be \defterm{Boolean} if for any quadruple $abcd\in X^4$,\newline $\Aline ab\cap\Aline cd=\varnothing=\Aline bc\cap\Aline ad$ implies $\Aline ac\cap\Aline bd=\varnothing$.
\end{definition}

\begin{theorem}[{\cite[6.8.8]{Banakh}}]\label{t:Boolean<=>} A $3$-ranked liner $X$ is Boolean if and only if every parallelogram in $X$ is Boolean.
\end{theorem}

\begin{example} Let $V$ be a vector space over a field $F$. The space $V$ endowed with its canonical line structure is a Boolean liner if and only fi the field $F$ has characteristic two. 
\end{example}

\section{Hyper-Bolyai liners}

Motivated by the problem of extending triangles to parallelo\-grams, we introduce and study in this section the class of hyper-Bolyai liners, intermediate between the Playfair and Bolyai liners. 

Let us recall that a liner $X$ is \defterm{Bolyai} if for any plane $P\subseteq X$, line $L\subseteq P$ and point $p\in P\setminus L$ there exists a line $\Lambda$ such that $p\in \Lambda\subseteq P\setminus L$.

\begin{definition}\label{d:hyper-Bolyai} A liner $X$ is defined to be \defterm{hyper-Bolyai} if 
for any plane $P\subseteq X$, concurrent lines $A,B\subseteq P$ and point $p\in P\setminus B$, there exists a line $\Lambda\subseteq X$ such that $p\in \Lambda$, $\Lambda\parallel B$, and $|\Lambda\cap A|=1$. 
\end{definition}

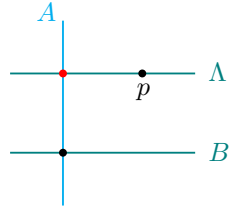
\begin{picture}(100,80)(-150,-25)
\linethickness{0.7pt}
\put(0,0){\color{teal}\line(1,0){70}}
\put(75,-3){\color{teal}$B$}
\put(0,30){\color{teal}\line(1,0){70}}
\put(75,27){\color{teal}$\Lambda$}
\put(20,-20){\color{cyan}\line(0,1){70}}
\put(10,50){\color{cyan}$A$}
\put(20,30){\color{red}\circle*{3}}
\put(50,30){\circle*{3}}
\put(48,22){$p$}
\put(20,0){\circle*{3}}
\end{picture}

Hyper-Bolyai liners relate to some other types of liners as follows.
$$
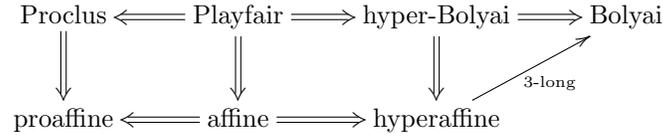
\xymatrix{
\mbox{Proclus}\ar@{=>}[d]&\mbox{Playfair}\ar@{=>}[l]\ar@{=>}[r]\ar@{=>}[d]&\mbox{hyper-Bolyai}\ar@{=>}[r]\ar@{=>}[d]&\mbox{Bolyai}\\
\mbox{proaffine}&\mbox{affine}\ar@{=>}[l]\ar@{=>}[r]&\mbox{hyperaffine}\ar_{\tiny\mbox{$3$-long}}[ru]
}
$$
Non-trivial implications in this diagram are proved in Theorem~\ref{t:Proclus<=>}, \ref{t:Playfair<=>}, and the following proposition.

\begin{proposition}[{\cite[6.9.2]{Banakh}}]\label{p:hyper-Bolyai-interplay}
\begin{enumerate}\setlength{\itemsep}{-0pt}\setlength{\parskip}{1pt}
\item[\textup{(1)}] Every hyper-Bolyai liner is Bolyai.
\item[\textup{(2)}] Every hyper-Bolyai liner is hyperaffine.
\item[\textup{(3)}] Every hyperaffine $3$-long is Bolyai.
\item[\textup{(4)}] A liner is Playfair if and only if it is Proclus and hyper-Bolyai.
\end{enumerate}
\end{proposition}

\begin{remark} By Example~\ref{ex:non-hyperaffine}, there exists a $4$-long Bolyai plane, which is not hyperaffine. By Example~\ref{ex:hyperaffine-nonPlayfair}, there exists a  non-Playfair hyperaffine $4$-long ranked plane, which is hyper-Bolyai.
\end{remark}

The following theorem characterizes hyper-Bolyai liners, in the spirit of the characterizations of Proclus and Playfair liners given in Theorems \ref{t:Proclus<=>} and \ref{t:Playfair<=>}.

\begin{theorem}[{\cite[6.6.5]{Banakh}}]\label{t:hyper-Bolyai<=>} A liner $X$ of cardinality $|X|>2$ is hyper-Bolyai if and only if $X$ is $3$-long, hyperaffine, and $3$-ranked.
\end{theorem}

The following proposition (which motivated our study of hyper-Bolyai liners) show that triangles in hyper-Bolyai liners can be extended to parallelograms.

\begin{proposition}[{\cite[6.9.6]{Banakh}}]\label{p:hyper-Bolyai=>parallelogram} For every triangle $abc$ in a hyper-Bolyai liner $X$, there exists a point $d\in X$ such that $abcd$ is a parallelogram.
\end{proposition}

\begin{remark} By Example~\ref{ex:non-hyperaffine}, there exists a $4$-parallel plane $X$ with $|X|_2=4$ and $|X|=|X|_3=25$, which is not hyperaffine and hence not hyper-Bolyai.  By Proposition~\ref{p:k-parallel=>parallelogram}, every triangle in $X$ can be completed to a parallelogram. This example shows that Proposition~\ref{p:hyper-Bolyai=>parallelogram} cannot be turned into a characterization of hyper-Bolyai liners.
\end{remark}

\end{document}